\documentclass[12pt]{article}

\usepackage{amsmath,amssymb,amsthm,amscd,a4wide,upref}

\usepackage[enableskew]{youngtab}

\usepackage{hyperref}
\hypersetup{colorlinks,citecolor=blue,filecolor=black,linkcolor=blue,urlcolor=blue}
\usepackage{enumerate}
\usepackage{mathtools}
\usepackage{enumitem}

\usepackage{lmodern}     
\usepackage[T1]{fontenc}

\begin{document}


\newcommand{\ad}{{\rm ad}}
\newcommand{\cri}{{\rm cri}}
\newcommand{\row}{{\rm row}}
\newcommand{\col}{{\rm col}}
\newcommand{\Ann}{{\rm{Ann}\ts}}
\newcommand{\End}{{\rm{End}\ts}}
\newcommand{\Rep}{{\rm{Rep}\ts}}
\newcommand{\Hom}{{\rm{Hom}}}
\newcommand{\Mat}{{\rm{Mat}}}
\newcommand{\ch}{{\rm{ch}\ts}}
\newcommand{\chara}{{\rm{char}\ts}}
\newcommand{\diag}{{\rm diag}}
\newcommand{\st}{{\rm st}}
\newcommand{\non}{\nonumber}
\newcommand{\wt}{\widetilde}
\newcommand{\wh}{\widehat}
\newcommand{\ol}{\overline}
\newcommand{\ot}{\otimes}
\newcommand{\la}{\lambda}
\newcommand{\La}{\Lambda}
\newcommand{\De}{\Delta}
\newcommand{\al}{\alpha}
\newcommand{\be}{\beta}
\newcommand{\ga}{\gamma}
\newcommand{\Ga}{\Gamma}
\newcommand{\ep}{\epsilon}
\newcommand{\ka}{\kappa}
\newcommand{\vk}{\varkappa}
\newcommand{\si}{\sigma}
\newcommand{\vs}{\varsigma}
\newcommand{\vp}{\varphi}
\newcommand{\ta}{\theta}
\newcommand{\de}{\delta}
\newcommand{\ze}{\zeta}
\newcommand{\om}{\omega}
\newcommand{\Om}{\Omega}
\newcommand{\ee}{\epsilon^{}}
\newcommand{\su}{s^{}}
\newcommand{\hra}{\hookrightarrow}
\newcommand{\ve}{\varepsilon}
\newcommand{\pr}{^{\tss\prime}}
\newcommand{\ts}{\,}
\newcommand{\vac}{\mathbf{1}}
\newcommand{\vacu}{|0\rangle}
\newcommand{\di}{\partial}
\newcommand{\qin}{q^{-1}}
\newcommand{\tss}{\hspace{1pt}}
\newcommand{\Sr}{ {\rm S}}
\newcommand{\U}{ {\rm U}}
\newcommand{\BL}{ {\overline L}}
\newcommand{\BE}{ {\overline E}}
\newcommand{\BP}{ {\overline P}}
\newcommand{\AAb}{\mathbb{A}\tss}
\newcommand{\CC}{\mathbb{C}\tss}
\newcommand{\KK}{\mathbb{K}\tss}
\newcommand{\QQ}{\mathbb{Q}\tss}
\newcommand{\SSb}{\mathbb{S}\tss}
\newcommand{\TT}{\mathbb{T}\tss}
\newcommand{\ZZ}{\mathbb{Z}\tss}
\newcommand{\DY}{ {\rm DY}}
\newcommand{\X}{ {\rm X}}
\newcommand{\Y}{ {\rm Y}}
\newcommand{\Z}{{\rm Z}}
\newcommand{\ZX}{{\rm ZX}}
\newcommand{\Ac}{\mathcal{A}}
\newcommand{\Lc}{\mathcal{L}}
\newcommand{\Mc}{\mathcal{M}}
\newcommand{\Pc}{\mathcal{P}}
\newcommand{\Qc}{\mathcal{Q}}
\newcommand{\Rc}{\mathcal{R}}
\newcommand{\Sc}{\mathcal{S}}
\newcommand{\Tc}{\mathcal{T}}
\newcommand{\Bc}{\mathcal{B}}
\newcommand{\Ec}{\mathcal{E}}
\newcommand{\Fc}{\mathcal{F}}
\newcommand{\Gc}{\mathcal{G}}
\newcommand{\Hc}{\mathcal{H}}
\newcommand{\Uc}{\mathcal{U}}
\newcommand{\Vc}{\mathcal{V}}
\newcommand{\Wc}{\mathcal{W}}
\newcommand{\Yc}{\mathcal{Y}}
\newcommand{\Cl}{\mathcal{C}l}
\newcommand{\Ar}{{\rm A}}
\newcommand{\Br}{{\rm B}}
\newcommand{\Ir}{{\rm I}}
\newcommand{\Fr}{{\rm F}}
\newcommand{\Jr}{{\rm J}}
\newcommand{\Or}{{\rm O}}
\newcommand{\GL}{{\rm GL}}
\newcommand{\Spr}{{\rm Sp}}
\newcommand{\Rr}{{\rm R}}
\newcommand{\Zr}{{\rm Z}}
\newcommand{\gl}{\mathfrak{gl}}
\newcommand{\middd}{{\rm mid}}
\newcommand{\ev}{{\rm ev}}
\newcommand{\Pf}{{\rm Pf}}
\newcommand{\Norm}{{\rm Norm\tss}}
\newcommand{\oa}{\mathfrak{o}}
\newcommand{\spa}{\mathfrak{sp}}
\newcommand{\osp}{\mathfrak{osp}}
\newcommand{\f}{\mathfrak{f}}
\newcommand{\g}{\mathfrak{g}}
\newcommand{\h}{\mathfrak h}
\newcommand{\n}{\mathfrak n}
\newcommand{\m}{\mathfrak m}
\newcommand{\z}{\mathfrak{z}}
\newcommand{\Zgot}{\mathfrak{Z}}
\newcommand{\p}{\mathfrak{p}}
\newcommand{\sll}{\mathfrak{sl}}
\newcommand{\agot}{\mathfrak{a}}
\newcommand{\bgot}{\mathfrak{b}}
\newcommand{\qdet}{ {\rm qdet}\ts}
\newcommand{\Ber}{ {\rm Ber}\ts}
\newcommand{\HC}{ {\mathcal HC}}
\newcommand{\cdet}{{\rm cdet}}
\newcommand{\rdet}{{\rm rdet}}
\newcommand{\tr}{ {\rm tr}}
\newcommand{\gr}{ {\rm gr}\ts}
\newcommand{\str}{ {\rm str}}
\newcommand{\loc}{{\rm loc}}
\newcommand{\Gr}{{\rm G}}
\newcommand{\sgn}{ {\rm sgn}\ts}
\newcommand{\sign}{{\rm sgn}}
\newcommand{\ba}{a}
\newcommand{\bb}{b}
\newcommand{\bi}{i}
\newcommand{\bj}{j}
\newcommand{\bk}{k}
\newcommand{\bl}{l}
\newcommand{\bp}{p}
\newcommand{\hb}{\mathbf{h}}
\newcommand{\Sym}{\mathfrak S}
\newcommand{\fand}{\quad\text{and}\quad}
\newcommand{\Fand}{\qquad\text{and}\qquad}
\newcommand{\For}{\qquad\text{or}\qquad}
\newcommand{\for}{\quad\text{or}\quad}
\newcommand{\grpr}{{\rm gr}^{\tss\prime}\ts}
\newcommand{\degpr}{{\rm deg}^{\tss\prime}\tss}
\newcommand{\bideg}{{\rm bideg}\ts}

\renewcommand{\theequation}{\arabic{section}.\arabic{equation}}

\numberwithin{equation}{section}

\newtheorem{thm}{Theorem}[section]
\newtheorem{lem}[thm]{Lemma}
\newtheorem{prop}[thm]{Proposition}
\newtheorem{cor}[thm]{Corollary}
\newtheorem{conj}[thm]{Conjecture}
\newtheorem*{mthm}{Main Theorem}
\newtheorem*{mthma}{Theorem A}
\newtheorem*{mthmb}{Theorem B}
\newtheorem*{mthmc}{Theorem C}
\newtheorem*{mthmd}{Theorem D}

\theoremstyle{definition}
\newtheorem{defin}[thm]{Definition}

\theoremstyle{remark}
\newtheorem{remark}[thm]{Remark}
\newtheorem{example}[thm]{Example}
\newtheorem{examples}[thm]{Examples}

\newcommand{\bth}{\begin{thm}}
\renewcommand{\eth}{\end{thm}}
\newcommand{\bpr}{\begin{prop}}
\newcommand{\epr}{\end{prop}}
\newcommand{\ble}{\begin{lem}}
\newcommand{\ele}{\end{lem}}
\newcommand{\bco}{\begin{cor}}
\newcommand{\eco}{\end{cor}}
\newcommand{\bde}{\begin{defin}}
\newcommand{\ede}{\end{defin}}
\newcommand{\bex}{\begin{example}}
\newcommand{\eex}{\end{example}}
\newcommand{\bes}{\begin{examples}}
\newcommand{\ees}{\end{examples}}
\newcommand{\bre}{\begin{remark}}
\newcommand{\ere}{\end{remark}}
\newcommand{\bcj}{\begin{conj}}
\newcommand{\ecj}{\end{conj}}

\newcommand{\bal}{\begin{aligned}}
\newcommand{\eal}{\end{aligned}}
\newcommand{\beq}{\begin{equation}}
\newcommand{\eeq}{\end{equation}}
\newcommand{\ben}{\begin{equation*}}
\newcommand{\een}{\end{equation*}}

\newcommand{\bpf}{\begin{proof}}
\newcommand{\epf}{\end{proof}}

\def\beql#1{\begin{equation}\label{#1}}

\newcommand{\Res}{\mathop{\mathrm{Res}}}

\title{\Large\bf Representations of the orthosymplectic Yangian}

\author{A. I. Molev}

\date{} 
\maketitle


\begin{abstract}
We give a complete description of the finite-dimensional irreducible representations of
the Yangian associated with the orthosymplectic Lie superalgebra $\osp_{1|2}$.
The representations are parameterized by monic polynomials in one variable, they are
classified in terms of highest weights. We give explicit constructions
of a family of elementary modules of the Yangian and show that a wide class
of irreducible representations can be produced by taking
tensor products of the elementary modules.
\end{abstract}



%

\section{Introduction}\label{sec:int}
\setcounter{equation}{0}

It is well-known that the classification of finite-dimensional irreducible representations of
the orthosymplectic Lie superalgebras $\osp_{M|2n}$ is more complicated than that of its general linear
counterpart $\gl_{M|N}$; see e.g. books \cite{cw:dr} and \cite{m:ls}. It is therefore
not surprising that
such a disparity extends to the respective Yangians.

The Yangian for $\gl_{M|N}$ was introduced by Nazarov~\cite{n:qb} and its finite-dimensional
irreducible representations were classified by Zhang~\cite{zh:sy}. The orthosymplectic Yangians
$\Y(\osp_{M|2n})$ were
introduced by Arnaudon {\it et al.\/}~\cite{aacfr:rp}, and the general
classification problem still remains open. Our goal in this paper is to describe
finite-dimensional irreducible representations of
the Yangian $\Y(\osp_{1|2})$.
One can expect that this simplest case will be instrumental
in solving the general classification problem. In particular,
an extension to $\osp_{1|2n}$ appears to be straightforward.

According to \cite{aacfr:rp}, the $R$-matrix originated
in \cite{zz:rf} leads to the $RTT$-type definition
of the extended Yangian $\X(\osp_{M|2n})$. The Yangian $\Y(\osp_{M|2n})$ can be regarded
as the quotient of the extended Yangian by the ideal generated by elements of the center, implying that
the finite-dimensional irreducible representations of these two algebras are essentially the same.
An explicit description of the center and the
Hopf algebra structure were given in \cite{aacfr:rp}.
In the subsequent work \cite{aacfr:sy}, a Drinfeld-type presentation of the Yangian $\Y(\osp_{1|2})$
was produced, the double Yangian was constructed and its universal $R$-matrix
was calculated in an explicit form. Applications of the orthosymplectic Yangians to spin chain models
were discussed in \cite{aacfr:ba}. More recently, linear and quadratic evaluations
in the Yangian $\Y(\osp_{M|2n})$
were investigated in \cite{fikk:yy} and \cite{ikk:yb}.

By analogy with the $R$-matrix form of its orthogonal and symplectic counterparts,
every finite-dimensional irreducible representation of the
orthosymplectic Yangian is a highest weight representation; cf. \cite{amr:rp}.
For the extended Yangian
$\X(\osp_{1|2})$, this means that
the representations are determined by the triples $(\la_1(u),\la_2(u),\la_3(u))$ of formal series
$\la_i(u)\in 1+u^{-1}\CC[[u^{-1}]]$ which should satisfy the consistency condition
\beql{ylarel}
\la_1(u)\la_{3}(u+1/2)=\la_{2}(u)\la_{2}(u+1/2).
\eeq
The key step in the classification is to find
the conditions on such triples for
the corresponding representation
to be finite-dimensional. Here we need a super-extension of the approach
originally used by Tarasov~\cite{t:im} and
which had already been adapted and applied to the twisted Yangians associated with the classical
Lie algebras; see \cite[Ch.~3 \& 4]{m:yc}. This leads to the following result.

\begin{mthm}\label{thm:yclassi}
The finite-dimensional irreducible representation of the algebra $\X(\osp_{1|2})$
associated with the highest weight $(\la_1(u),\la_2(u),\la_3(u))$
is finite-dimensional if and only if
\beql{ydomire}
\frac{\la_2(u)}{\la_{1}(u)}=\frac{P(u+1)}{P(u)}
\eeq
for some monic polynomial $P(u)$ in $u$. Hence, the finite-dimensional irreducible
representations of the Yangian $\Y(\osp_{1|2})$ are parameterized by
monic polynomials $P(u)$.
\end{mthm}

The parametrization turns out to be the same as for the representations of the Yangian $\Y(\sll_2)$
given in \cite{d:nr}, and we will call $P(u)$ the {\em Drinfeld polynomial} of the representation.
The appearance of the consistency conditions \eqref{ylarel}
is best explained by the Gauss decomposition
of the generator matrix for $\X(\osp_{1|2})$ and can be compared
with the similar conditions for the Yangian
$\Y(\oa_3)$; see \cite{jl:ib}, \cite[Sec.~5.3]{jlm:ib}, which were
derived in a different way in \cite{amr:rp}. However, unlike the case of
the Yangian $\Y(\sll_2)$ (or $\Y(\oa_3)$),
there is no epimorphism $\X(\osp_{1|2})\to \U(\osp_{1|2})$, so that in general, representations of
$\osp_{1|2}$
do not extend to the Yangian.

An essential step in the proof of the Main Theorem,
which will be completed in Sec.~\ref{subsec:tp},
 is the analysis of the {\em elementary modules}
$L(\al,\be)$ over $\X(\osp_{1|2})$ associated with the highest weights of the form
\beql{hwelem}
\la_1(u)=\frac{u+\al}{u+\be},\qquad \la_2(u)=1,\qquad \la_3(u)=\frac{u+\be-1/2}{u+\al-1/2}.
\eeq
The corresponding {\em small Verma module} $M(\al,\be)$ turns out to be irreducible
if and only if $\be-\al$ and $\be-\al+1/2$ are not nonnegative integers.
The elementary modules $L(\al,\be)$
are the irreducible quotients of $M(\al,\be)$ and so they split into three families,
according to these conditions.
The module $L(\al,\be)$ is finite-dimensional if and only if $\be-\al\in\ZZ_+$.
In this case, when regarded as an $\osp_{1|2}$-module, $L(\al,\be)$
decomposes into the direct sum
\ben
L(\al,\be)\cong\bigoplus_{p=0}^{\lfloor\frac{\be-\al}{2}\rfloor} V(\be-\al-2\tss p),
\een
where $V(\mu)$ denotes the $2\mu+1$-dimensional $\osp_{1|2}$-module
with the highest weight $\mu\in\ZZ_+$.
In particular,
\ben
\dim L(\al,\be)=\binom{\be-\al+2}{2}.
\een

We construct a basis of each small Verma module $M(\al,\be)$ and give explicit formulas for
the action of the generators of $\X(\osp_{1|2})$. This leads to a corresponding description of all
elementary modules.
We show that, up to twisting by a multiplication automorphism of $\X(\osp_{1|2})$, every
finite-dimensional irreducible representation of this algebra is isomorphic to the quotient
of the submodule of the tensor product module of the form
\beql{tenprele}
L(\al_1,\be_1)\ot\dots\ot L(\al_k,\be_k),
\eeq
generated by the tensor product of the highest vectors. The final step is to investigate
irreducibility conditions for such tensor products.

In the case of the Yangian $\Y(\sll_2)$, an irreducibility criterion for tensor products of
evaluation modules
was given by Chari and Pressley~\cite{cp:yr}; see also \cite[Ch.~3]{m:yc}.
Such tensor products exhaust all finite-dimensional irreducible $\Y(\sll_2)$-modules.
This property turns out not to extend to representations of the Yangian for $\osp_{1|2}$;
see Example~\ref{ex:tpr} below. A wide class of irreducible modules can still be constructed
explicitly via tensor products of the form \eqref{tenprele}; see Theorem~\ref{thm:suffc}.

\section{Sign conventions, definitions and preliminaries}
\label{sec:ns}

Consider the three-dimensional space $\CC^{1|2}$ over $\CC$ with basis elements $e_1,e_2,e_3$
where we assume a $\ZZ_2$-gradation defined by setting that the vector $e_2$ is even and
the vectors $e_1$ and $e_3$ are odd. It will be convenient to use the involution
on the set $\{1,2,3\}$ defined by $i\mapsto i'=4-i$.
The endomorphism algebra $\End\CC^{1|2}$ then gets a $\ZZ_2$-gradation with
the parity of the matrix unit $e_{ij}$ found by
$\bi+\bj\mod 2$.

We will consider $3\times 3$ matrices with entries in superalgebras. All such
matrices will be even so that the $(i,j)$ entry will have the parity $\bi+\bj\mod 2$.
We want to posit that the product of two matrices with entries in a superalgebra $\Ac$
is calculated in the standard way (without any additional signs). Accordingly, the algebra of
even matrices over $\Ac$ will be identified with the tensor product algebra
$\End\CC^{1|2}\ot\Ac$. A matrix $A=[a_{ij}]$ will be regarded as the element
\ben
A=\sum_{i,j=1}^{3}e_{ij}\ot a_{ij}(-1)^{\bi\tss\bj+\bj}\in \End\CC^{1|2}\ot\Ac.
\een

We will use the involutive matrix {\em super-transposition} $t$ defined by
$(A^t)_{ij}=A_{j'i'}(-1)^{\bi\bj+\bj}\tss\ta_i\ta_j$,
where we set
\ben
\ta_1=\ta_2=1\Fand \ta_3=-1.
\een
This super-transposition is associated with the bilinear form on the space $\CC^{1|2}$
defined by the anti-diagonal matrix $G=[\de_{ij'}\tss\ta_i]$.
We will also regard $t$ as the linear map
\beql{suptra}
t:\End\CC^{1|2}\to \End\CC^{1|2}, \qquad
e_{ij}\mapsto e_{j'i'}(-1)^{\bi\bj+\bi}\tss\ta_i\ta_j.
\eeq
In the case of multiple tensor products of the endomorphism algebras,
we will indicate by $t_a$ the map \eqref{suptra}
acting on the $a$-th copy of $\End\CC^{1|2}$.

A standard basis of the general linear Lie superalgebra $\gl_{1|2}$ is formed by elements $E_{ij}$
of the parity $\bi+\bj\mod 2$ for $1\leqslant i,j\leqslant 3$ with the commutation relations
\ben
[E_{ij},E_{kl}]
=\de_{kj}\ts E_{i\tss l}-\de_{i\tss l}\ts E_{kj}(-1)^{(\bi+\bj)(\bk+\bl)}.
\een
We will regard the orthosymplectic Lie superalgebra $\osp_{1|2}$
associated with the bilinear from defined by $G$ as the subalgebra
of $\gl_{1|2}$ spanned by the elements
\ben
F_{ij}=E_{ij}-E_{j'i'}(-1)^{\bi\tss\bj+\bi}\ts\ta_i\ta_j.
\een

For any given $\mu\in\CC$ we will denote by $V(\mu)$ the irreducible highest weight module over
$\osp_{1|2}$ generated by a nonzero vector $\xi$ such that
\ben
F_{11}\xi=\mu\ts\xi\Fand F_{12}\tss \xi=0.
\een
The module $V(\mu)$ is finite-dimensional if and only if $\mu\in\ZZ_+$. In that case,
$\dim V(\mu)=2\mu+1$.

Introduce the permutation operator $P$ as the element
\ben
P=\sum_{i,j=1}^3 e_{ij}\ot e_{ji}(-1)^{\bj}\in \End\CC^{1|2}\ot\End\CC^{1|2}
\een
and set
\ben
Q=P^{\tss t_1}=P^{\tss t_2}=\sum_{i,j=1}^3 e_{ij}\ot e_{i'j'}(-1)^{\bi\bj}\ts\ta_i\ta_j
\in \End\CC^{1|2}\ot\End\CC^{1|2}.
\een
The rational function in $u$ given by
\ben
R(u)=1-\frac{P}{u}+\frac{Q}{u-\ka},\qquad \ka=-3/2,
\een
is an $R$-{\em matrix}, it satisfies
the {\em Yang--Baxter equation} as originally found in \cite{zz:rf}. The $R$-matrices produced
in that paper are known to extend to the Brauer algebra so that
the Yang--Baxter equation can be verified
by taking a suitable Brauer algebra representation
in tensor products of the $\ZZ_2$-graded spaces; cf. \cite{fikk:yy}, \cite{imo:nf}.

Following \cite{aacfr:rp},
define the {\it extended Yangian\/}
$\X(\osp_{1|2})$
as an associative superalgebra with generators
$t_{ij}^{(r)}$ of parity $\bi+\bj\mod 2$, where $1\leqslant i,j\leqslant 3$ and $r=1,2,\dots$,
satisfying certain quadratic relations. In order to write them down,
introduce the formal series
\beql{tiju}
t_{ij}(u)=\de_{ij}+\sum_{r=1}^{\infty}t_{ij}^{(r)}\ts u^{-r}
\in\X(\osp_{1|2})[[u^{-1}]]
\eeq
and combine them into the $3\times 3$ matrix $T(u)=[t_{ij}(u)]$ so that
\ben
T(u)=\sum_{i,j=1}^3 e_{ij}\ot t_{ij}(u)(-1)^{\bi\tss\bj+\bj}
\in \End\CC^{1|2}\ot \X(\osp_{1|2})[[u^{-1}]].
\een
Consider the algebra
$\End\CC^{1|2}\ot\End\CC^{1|2}\ot \X(\osp_{1|2})[[u^{-1}]]$
and introduce its elements $T_1(u)$ and $T_2(u)$ by
\ben
T_1(u)=\sum_{i,j=1}^3 e_{ij}\ot 1\ot t_{ij}(u)(-1)^{\bi\tss\bj+\bj},\qquad
T_2(u)=\sum_{i,j=1}^3 1\ot e_{ij}\ot t_{ij}(u)(-1)^{\bi\tss\bj+\bj}.
\een
The defining relations for the superalgebra $\X(\osp_{1|2})$ can then
be written in the form of the $RTT$-{\em relation}
\beql{RTT}
R(u-v)\ts T_1(u)\ts T_2(v)=T_2(v)\ts T_1(u)\ts R(u-v).
\eeq
As shown in \cite{aacfr:rp}, the product $T(u)\ts T^{\tss t}(u-\ka)$ is a scalar matrix with
\beql{ttra}
T(u-\ka)\ts T^{\tss t}(u)=c(u)1,
\eeq
where $c(u)$ is a series in $u^{-1}$. All its coefficients belong to
the center $\ZX(\osp_{1|2})$ of $\X(\osp_{1|2})$ and generate the center.

The {\em Yangian} $\Y(\osp_{1|2})$
is defined as the subalgebra of
$\X(\osp_{1|2})$ which
consists of the elements stable under
the automorphisms
\beql{muf}
t_{ij}(u)\mapsto f(u)\ts t_{ij}(u)
\eeq
for all series
$f(u)\in 1+u^{-1}\CC[[u^{-1}]]$.
We have the tensor product decomposition
\beql{tensordecom}
\X(\osp_{1|2})=\ZX(\osp_{1|2})\ot \Y(\osp_{1|2}).
\eeq
The Yangian $\Y(\osp_{1|2})$ can be equivalently defined as the quotient
of $\X(\osp_{1|2})$
by the relation
\ben
T(u-\ka)\ts T^{\tss t}(u)=1.
\een

We will also use a more explicit form of the defining relations \eqref{RTT}
written in terms of the series \eqref{tiju} as follows:
\begin{align}
[\tss t_{ij}(u),t_{kl}(v)]&=\frac{1}{u-v}
\big(t_{kj}(u)\ts t_{il}(v)-t_{kj}(v)\ts t_{il}(u)\big)
(-1)^{\bi\tss\bj+\bi\tss\bk+\bj\tss\bk}
\non\\
{}&-\frac{1}{u-v-\kappa}
\Big(\de_{k i\pr}\sum_{p=1}^3\ts t_{pj}(u)\ts t_{p'l}(v)
(-1)^{\bi+\bi\tss\bj+\bj\tss\bp}\ts\ta_i\ta_p
\label{defrel}\\
&\qquad\qquad\qquad
{}-\de_{l j\pr}\sum_{p=1}^3\ts t_{k\tss p'}(v)\ts t_{ip}(u)
(-1)^{\bj+\bp+\bi\tss\bk+\bj\tss\bk+\bi\tss\bp}\ts\ta_j\ta_p\Big).
\non
\end{align}
The coefficients of the series $t_{11}(u), t_{12}(u), t_{21}(u)$ and $c(u)$
generate the algebra $\X(\osp_{1|2})$.
The mapping
$
t_{ij}(u)\mapsto t_{ij}(-u)
$
defines an anti-automorphism of $\X(\osp_{1|2})$, while each of the mappings
\beql{shift}
t_{ij}(u)\mapsto t_{ij}(u+a),\quad a\in \CC,
\eeq
and
$
t_{ij}(u)\mapsto t_{i'j'}(u)\ts\ta_i\ta_j
$
defines an automorphism. Consider their composition to define the anti-automorphism
\beql{om}
\om:t_{ij}(u)\mapsto t_{i'j'}(-u+1/2)\ts\ta_i\ta_j.
\eeq

The universal enveloping algebra $\U(\osp_{1|2})$ can be regarded as a subalgebra of
$\X(\osp_{1|2})$ via the embedding
\beql{emb}
F_{ij}\mapsto \frac12\big(t_{ij}^{(1)}-t_{j'i'}^{(1)}(-1)^{\bj+\bi\bj}\ts\ta_i\ta_j\big)(-1)^{\bi}.
\eeq
This fact relies on the Poincar\'e--Birkhoff--Witt theorem for the orthosymplectic Yangian
which was essentially pointed out in \cite{aacfr:rp} and \cite{aacfr:sy}, as
the associated graded algebra
for $\Y(\osp_{1|2})$ is isomorphic to $\U(\osp_{1|2}[u])$.
It states that given any total ordering on the set of generators
$t_{ij}^{(r)}$ with $i+j\leqslant 4$ and $r\geqslant 1$,
the ordered monomials
in the generators with the powers of odd generators not exceeding $1$,
form a basis of $\X(\osp_{1|2})$.
A detailed proof can be carried over by
adapting the arguments of \cite[Sec.~3]{amr:rp} to the super case in a straightforward way
with the use of the vector representation recalled below in \eqref{vectre}.

The extended Yangian $\X(\osp_{1|2})$ is a Hopf algebra with the coproduct
defined by
\beql{Delta}
\De: t_{ij}(u)\mapsto \sum_{k=1}^3 t_{ik}(u)\ot t_{kj}(u).
\eeq
For the image of the series $c(u)$ we have $\De:c(u)\mapsto c(u)\ot c(u)$ and so the Yangian
$\Y(\osp_{1|2})$ inherits the Hopf algebra structure from $\X(\osp_{1|2})$.

\section{Gaussian generators}
\label{sec:gd}

A Drinfeld-type presentation of the Yangian for $\osp_{1|2}$ was given in \cite{aacfr:sy}
with the use of the Gauss decomposition of the matrix $T(u)$. We will be using
some calculations
produced therein and derive consistency relations for the Gaussian generators.

Apply the Gauss decomposition
to the generator matrix $T(u)$ for $\X(\osp_{1|2})$,
\beql{gd}
T(u)=F(u)\ts H(u)\ts E(u),
\eeq
where $F(u)$, $H(u)$ and $E(u)$ are uniquely determined matrices of the form
\ben
F(u)=\begin{bmatrix}
1&0&0\ts\\
f_{21}(u)&1&0\\
f_{31}(u)&f_{32}(u)&1
\end{bmatrix},
\qquad
E(u)=\begin{bmatrix}
\ts1&e_{12}(u)&e_{13}(u)\ts\\
\ts0&1&e_{23}(u)\\
0&0&1
\end{bmatrix},
\een
and $H(u)=\diag\ts\big[h_1(u),h_2(u),h_3(u)\big]$.
Explicit formulas for the entries
of the matrices $F(u)$, $H(u)$ and $E(u)$ can be written with the use of the
Gelfand--Retakh quasideterminants \cite{gr:dm}; cf. \cite[Sec.~4]{jlm:ib}.
In particular, we have
\ben
h_1(u)=t_{11}(u),\qquad h_2(u)=\begin{vmatrix}
t_{11}(u)&t_{12}(u)\\
t_{21}(u)&\boxed{t_{22}(u)}
\end{vmatrix},\qquad
h_3(u)=\begin{vmatrix} t_{11}(u)&t_{12}(u)&t_{13}(u)\\[0.2em]
                         t_{21}(u)&t_{22}(u)&t_{23}(u)\\
                         t_{31}(u)&t_{32}(u)&\boxed{t_{33}(u)}
           \end{vmatrix},
\een
whereas
\ben
e_{12}(u)=h_1(u)^{-1}\ts t_{12}(u),\qquad
           e_{23}(u)=h_2(u)^{-1}\ts\begin{vmatrix}
t_{11}(u)&t_{13}(u)\\
t_{21}(u)&\boxed{t_{23}(u)}
\end{vmatrix},
\een
and
\ben
f_{21}(u)=t_{21}(u)\ts h_1(u)^{-1},\qquad
           f_{32}(u)=\begin{vmatrix}
t_{11}(u)&t_{12}(u)\\
t_{31}(u)&\boxed{t_{32}(u)}
\end{vmatrix}\ts h_2(u)^{-1}.
\een

\bpr\label{prop:gaussrel}
The following relations for the Gaussian generators hold:
\beql{ef}
e_{12}(u)=-e_{23}(u+1/2),\qquad f_{21}(u)=f_{32}(u+1/2),
\eeq
and
\beql{hoht}
h_1(u)\tss h_3(u+1/2)=h_2(u)\tss h_2(u+1/2).
\eeq
Moreover,
\beql{cu}
c(u)=h_1(u)\tss h_1(u+1)^{-1}\tss h_2(u+1)\tss h_2(u+3/2).
\eeq
\epr

\bpf
The argument is quite similar to the proof of the corresponding relations for the Gaussian generators
of $\Y(\oa_3)$ given in \cite{jl:ib}; see also \cite[Sec.~5.3]{jlm:ib}. We will outline a few key steps.

By inverting the matrices on both sides of \eqref{gd}, we get
\ben
T(u)^{-1}=E(u)^{-1}\tss H(u)^{-1}\tss F(u)^{-1}.
\een
On the other hand, relation \eqref{ttra} implies $T^{\tss t}(u)=c(u)\tss T(u-\ka)^{-1}$.
Hence, by equating the $(i,j)$ entries with $i,j=2,3$ in this matrix relation, we derive
\begin{align}
h_1(u)&=c(u)\tss h_3(u-\ka)^{-1},
\non\\[0.2em]
h_1(u)\tss e_{12}(u)&=-c(u)\tss e_{23}(u-\ka)\tss h_3(u-\ka)^{-1},
\label{hef}\\[0.3em]
f_{21}(u)\tss h_1(u)&=c(u)\tss h_3(u-\ka)^{-1}\tss f_{32}(u-\ka),
\non
\end{align}
and
\begin{align}
h_2(u)+f_{21}(u)\tss & h_1(u)\tss e_{12}(u)
\non\\
{}&=c(u)\big(h_2(u-\ka)^{-1}+e_{23}(u-\ka)\tss h_3(u-\ka)^{-1}\tss f_{32}(u-\ka)\big).
\label{twot}
\end{align}
Calculating as in \cite{aacfr:sy} and \cite{jl:ib}, we verify that the coefficients of the series
$h_1(u), h_2(u)$ and $h_3(u)$ pairwise commute. Furthermore,
we get
\ben
h_1(u)\tss e_{12}(u)=e_{12}(u+1)\tss h_1(u)\Fand h_1(u)\tss f_{21}(u+1)=f_{21}(u)\tss h_1(u)
\een
which together with relations \eqref{hef} imply the first two desired identities, where
we replaced $\ka$ by its value $-3/2$. They imply that
relation \eqref{twot} can be written in the form
\beql{hsi}
h_2(u)-c(u)\tss h_2(u-\ka)^{-1}=-\big[e_{12}(u+1),f_{21}(u)\big]\tss h_1(u).
\eeq
As a final step, use one more relation between the Gaussian generators,
\ben
\big[e_{12}(u),f_{21}(v)\big]=\frac{h_1(u)^{-1}\tss h_2(u)-h_1(v)^{-1}\tss h_2(v)}{u-v},
\een
so that eliminating $c(u)$ from \eqref{hsi} we come to \eqref{hoht}.
Relation \eqref{cu} follows by eliminating $h_3(u)$ from the first relation in \eqref{hef}
with the use of \eqref{hoht}.
\epf

Observe that the coefficients of the series $e_{12}(u)$ and $f_{21}(u)$ are stable under
all automorphisms \eqref{muf} and so belong to the subalgebra $\Y(\osp_{1|2})$
of $\X(\osp_{1|2})$.
Together with the coefficients of the series $h(u)=h_1(u)^{-1}h_2(u)$
they generate the Yangian $\Y(\osp_{1|2})$, and the defining relations for these generators
are given in \cite{aacfr:sy} in a slightly different setting.

\section{Yangian representations}
\label{sec:re}

\subsection{Highest weight representations}
\label{subsec:hw}

A representation $V$ of the algebra $\X(\osp_{1|2})$
is called a {\em highest weight representation}
if there exists a nonzero vector
$\xi\in V$ such that $V$ is generated by $\xi$,
\begin{alignat}{2}
t_{ij}(u)\ts\xi&=0 \qquad &&\text{for}
\quad 1\leqslant i<j\leqslant 3, \qquad \text{and}\non\\
t_{ii}(u)\ts\xi&=\la_i(u)\ts\xi \qquad &&\text{for}
\quad i=1,2,3,
\label{trianb}
\end{alignat}
for some formal series
\beql{laiu}
\la_i(u)\in 1+u^{-1}\CC[[u^{-1}]].
\eeq
The vector $\xi$ is called the {\em highest vector}
of $V$ and the triple $\la(u)=(\la_{1}(u),\la_2(u),\la_3(u))$
is the {\em highest weight\/} of $V$.

The quasideterminant formulas for the Gaussian generators $h_i(u)$ given in Sec.~\ref{sec:gd}
imply that the conditions \eqref{trianb} in the above definition can be replaced with
\ben
h_{i}(u)\ts\xi=\la_i(u)\ts\xi \qquad \text{for}
\quad i=1,2,3.
\een
Hence, relation \eqref{hoht} of
Proposition~\ref{prop:gaussrel} implies the consistency condition \eqref{ylarel}
for the components $\la_i(u)$ of the highest weight.

Given any triple $\la(u)=(\la_{1}(u),\la_2(u),\la_3(u))$
of formal series of the form \eqref{laiu} safisfying the consistency
condition \eqref{ylarel}, define
the {\em Verma module} $M(\la(u))$ as the quotient of $\X(\osp_{1|2})$ by
the left ideal generated by all coefficients of the series $t_{ij}(u)$
with $1\leqslant i<j\leqslant 3$, and $t_{ii}(u)-\la_i(u)$ for
$i=1,2,3$. The Poincar\'e--Birkhoff--Witt theorem implies that the Verma module $M(\la(u))$
is nonzero and we denote by $L(\la(u))$ its irreducible quotient.

\bre\label{rem:dp}
In terms of the Drinfeld presentation
of the Yangian $\Y(\osp_{1|2})$ given in \cite{aacfr:sy}, the highest vector conditions
take the form
\ben
e_{12}(u)\ts\xi=0\Fand h(u)\ts\xi=\mu(u)\ts\xi
\een
for a series $\mu(u)\in 1+u^{-1}\CC[[u^{-1}]]$.
The Main Theorem implies that the irreducible highest weight representation of $\Y(\osp_{1|2})$
associated with $\mu(u)$ is finite-dimensional if and only if
\ben
\mu(u)=\frac{P(u+1)}{P(u)}
\een
for some monic polynomial $P(u)$ in $u$.
\qed
\ere

\bpr\label{prop:fdhw}
Every finite-dimensional irreducible representation of the algebra $\X(\osp_{1|2})$
is isomorphic to $L(\la(u))$ for a certain highest weight $\la(u)$ satisfying \eqref{ylarel}.
\epr

\bpf
The argument is a straightforward adaptation of the proof of \cite[Thm~5.1]{amr:rp} to the super case,
which amounts to taking care of additional signs in the calculations.
\epf

To prove the Main Theorem, we only
need to determine which of the modules $L(\la(u))$ are finite-dimensional.
The first step is to establish some necessary conditions.

\bpr\label{prop:rat}
If the module $L(\la(u))$ is finite-dimensional, then
\ben
\frac{\la_1(u)}{\la_2(u)}=\frac{(u+\al_1)\dots (u+\al_k)}{(u+\be_1)\dots (u+\be_k)}
\een
for $k\in\ZZ_+$ and certain complex numbers $\al_i,\be_i$.
\epr

\bpf
We follow the proof of a similar property for the Yangian $\Y(\gl_2)$; see \cite[Prop.~3.3.1]{m:yc}.
Note that due to \eqref{ylarel}, the series $\la_3(u)$ is uniquely
determined by $\la_{1}(u)$ and $\la_2(u)$, and so we can parameterize
the highest weights by arbitrary pairs of series $\la(u)=(\la_{1}(u),\la_2(u))$,
omitting $\la_3(u)$.
By twisting the action of the extended Yangian $\X(\osp_{1|2})$ on $L(\la(u))$
by the automorphism \eqref{muf}
with $f(u)=\la_2(u)^{-1}$,
we get an $\X(\osp_{1|2})$-module isomorphic to
$L(\mu(u),1)$ with $\mu(u)=\la_1(u)/\la_2(u)$.
Let $\xi$ denote the highest vector of $L(\mu(u),1)$.
Since this representation is finite-dimensional,
the vectors $t_{21}^{(i)}\xi\in L(\mu(u),1)$ with $i\geqslant 1$ are
linearly dependent,
\ben
\sum_{i=1}^m c_i\ts t_{21}^{(i)}\xi=0
\een
with $c_i\in\CC$, assuming $c_m\ne 0$. Apply the operators
$t_{12}^{(r)}$ for all $r\geqslant 1$ to the linear combination on the left hand side
and take the coefficient of $\xi$. Since $t_{12}(u)\tss\xi=0$,
we get from the defining relations \eqref{defrel} that
\ben
t_{12}(u)\tss t_{21}(v)\tss\xi=\frac{1}{u-v}
\big(t_{22}(u)\ts t_{11}(v)-t_{22}(v)\ts t_{11}(u)\big)\tss\xi=-\frac{\mu(u)-\mu(v)}{u-v}\tss\xi.
\een
Hence, writing
\ben
\mu(u)=1+\mu^{(1)}u^{-1}+\mu^{(2)}u^{-2}+\dots,\qquad \mu^{(i)}\in\CC,
\een
we derive
\ben
t_{12}^{(r)}\ts t_{21}^{(i)}\xi=\mu^{(r+i-1)}\tss \xi.
\een
Therefore, for all $r\geqslant 1$ we have the relations
\ben
\sum_{i=1}^m c_i\tss \mu^{(r+i-1)}=0.
\een
They imply
\ben
\mu(u)\tss(c_1+c_2\tss u+\dots+c_m\tss u^{m-1})
=(b_1+b_2\tss u+\dots+b_m\tss u^{m-1})
\een
so that $\mu(u)$ can be written as a rational function in $u$, as required.
\epf

We will use the name {\em elementary module} for the module $L(\la(u))$ with
\beql{hwe}
\la_1(u)=\frac{u+\al}{u+\be}\Fand \la_2(u)=1
\eeq
and denote it by $L(\al,\be)$.

The Hopf algebra structure on the extended Yangian $\X(\osp_{1|2})$ allows us to regard
tensor products of the form
\beql{tenprelem}
L=L(\al_1,\be_1)\ot\dots\ot L(\al_k,\be_k)
\eeq
as $\X(\osp_{1|2})$-modules. Let $\xi^{(i)}$ denote the highest vector of $L(\al_i,\be_i)$.

\bpr\label{prop:subq}
The $\X(\osp_{1|2})$-module $L(\la(u))$ with
\beql{hwgen}
\la_1(u)=\frac{(u+\al_1)\dots (u+\al_k)}{(u+\be_1)\dots (u+\be_k)}\Fand \la_2(u)=1
\eeq
is isomorphic to the irreducible quotient of the submodule of $L$,
generated by
the tensor product of the highest vectors $\xi^{(1)}\ot\dots\ot \xi^{(k)}$.
\epr

\bpf
The coproduct formula \eqref{Delta}
implies that
the cyclic span $\X(\osp_{1|2})(\xi^{(1)}\ot\dots\ot \xi^{(k)})$
is a highest weight module
with the highest weight $(\la_1(u),\la_2(u))$ which implies the claim.
\epf

We will need to find the conditions for the elementary modules to be finite-dimensional and establish
some sufficient conditions for the module $L$ in \eqref{tenprelem} to be irreducible.

\subsection{Small Verma modules}
\label{subsec:rv}

Note that by twisting the action of the extended Yangian
on a highest weight module
with the highest weight \eqref{hwe} by the shift automorphism \eqref{shift} with $a=-\be$,
we get the corresponding module whose highest weight is found by shifting $\al\mapsto\al-\be$
and $\be\mapsto 0$. We will now assume that $\be=0$.
Let $\al\in \CC$ and consider the Verma module $M(\la(u))$
with
\beql{hwal}
\la_1(u)=\frac{u+\al}{u},\qquad\la_2(u)=1,\qquad \la_3(u)=\frac{u-1/2}{u+\al-1/2}.
\eeq
Let $K$
be the submodule of $M(\la(u))$ generated by all vectors of the form
\beql{gensub}
t_{21}^{(r)}\tss\xi\quad\text{for}\quad r\geqslant 2\Fand
\big(t_{31}^{(r)}+(\al-1/2)\tss t_{31}^{(r-1)}\big)\tss\xi\quad\text{for}\quad r\geqslant 3,
\eeq
where $\xi$ denotes the highest vector of the Verma module.
Introduce the {\em small Verma module} $M(\al)$ as the quotient $M(\la(u))/K$.
We will keep the notation $\xi$ for the image of the highest vector
of the Verma module in the quotient.
More general small Verma modules of the form $M(\al,\be)$ corresponding to the
highest weights \eqref{hwe} are then obtained by twisting
the modules $M(\al)$ by suitable automorphisms \eqref{shift}.

\bpr\label{prop:span}
The module $M(\al)$ is spanned by the vectors
\beql{vects}
t_{31}^{(2)\tss r}\tss t_{21}^{(1)\tss s}\ts \xi,\qquad r,s\in\ZZ_+.
\eeq
\epr

\bpf
By the Poincar\'e--Birkhoff--Witt theorem for the extended Yangian, the Verma module
$M(\la(u))$ has the basis
\beql{pbwb}
t_{31}^{(k_1)}\dots t_{31}^{(k_p)}\tss t_{21}^{(l_1)}\dots t_{21}^{(l_q)}\tss\xi,
\eeq
where $k_1\geqslant\dots\geqslant k_p\geqslant 1$ and $l_1>\dots>l_q\geqslant 1$. Hence, the induction
on the length of the monomial in \eqref{pbwb} reduces the argument to the verification
of the property that
the span of the vectors \eqref{vects} is stable under the action of the generators
$t_{31}^{(k)}$ and $t_{21}^{(l)}$.

The defining relations \eqref{defrel} imply that $[t_{31}^{(k)},t_{31}^{(m)}]=0$
and $[t_{31}^{(k)},t_{21}^{(1)}]=0$ for all $k,m$. Therefore, for $k\geqslant 2$ in $M(\la(u))$
we have
\ben
t_{31}^{(k)}\tss t_{31}^{(2)\tss r}\tss t_{21}^{(1)\tss s}\ts \xi\equiv(-\al+1/2)^{k-2}\ts
t_{31}^{(2)\tss r+1}\tss t_{21}^{(1)\tss s}\ts \xi\mod K.
\een
The property is also clear for $k=1$ because $t_{31}^{(1)}=2\tss t_{21}^{(1)\tss 2}$.
Furthermore, since
\ben
[t_{21}^{(l)},t_{31}^{(2)}]=t_{21}^{(1)}\tss t_{31}^{(l)}-t_{31}^{(1)}\tss t_{21}^{(l)}
\een
and $[t_{21}^{(l)},t_{21}^{(1)}]=t_{31}^{(l)}$, the property for the generators
$t_{21}^{(l)}$ easily follows too.
\epf

We will regard $M(\al)$ as an $\osp_{1|2}$-module via the embedding \eqref{emb}.
We get the weight space decomposition
\ben
M(\al)=\bigoplus_{p=0}^{\infty} M(\al)_{-\al-p},
\een
where we define the weight subspaces of an arbitrary $\osp_{1|2}$-module $V$ by
\beql{wei}
V_{\ga}=\{v\in V\ |\ F_{11}\tss v=\ga\tss v\}.
\eeq
Proposition~\ref{prop:span} implies that
\beql{weineq}
\dim M(\al)_{-\al-p}\leqslant \lfloor p/2\rfloor+1.
\eeq

For all values $i,j\in\{1,2,3\}$
set $T_{ij}(u)=u\tss (u+\al-1/2)\ts t_{ij}(u)$. We will regard the coefficients
of these Laurent series in $u$ as operators in $M(\al)$.

\bpr\label{prop:pol}
All operators $T_{ij}(u)$ on the small Verma module $M(\al)$ are polynomials in $u$.
\epr

\bpf
Calculating modulo $K$, we get
\ben
t_{21}(u)\tss\xi=u^{-1}\tss t_{21}^{(1)}\tss \xi\Fand
t_{31}(u)\tss\xi=\big(u^{-1}\tss t_{31}^{(1)}+\frac{1}{u(u+\al-1/2)}t_{31}^{(2)}\big)\tss \xi
\een
so that the claim holds for the action of the operators $T_{21}(u)$ and $T_{31}(u)$ on $\xi$.
By acting on the vectors \eqref{vects} of the spanning set, we
note that the operator $T_{31}(u)$ commutes with $t_{31}^{(2)}$ and $t_{21}^{(1)}$, while
for the operator $T_{21}(u)$ we have the relations
\ben
[T_{21}(u),t_{31}^{(2)}]=t_{21}^{(1)}\tss T_{21}(u)-t_{31}^{(1)}\tss T_{31}(u)
\Fand [T_{21}(u),t_{21}^{(1)}]=T_{31}(u).
\een
Hence the property for the operators $T_{21}(u)$ and $T_{31}(u)$
follows by an obvious induction.

As a next step, consider the relations for the series $T_{11}(u)$ implied by
\eqref{defrel}:
\ben
[T_{11}(u),t_{21}^{(1)}]=T_{21}(u),\qquad [T_{11}(u),t_{31}^{(1)}]=2\tss T_{31}(u)
\een
and
\ben
[T_{11}(u),t_{31}^{(2)}]=T_{31}(u)\big(2\tss u+1/2+t_{11}^{(1)}\big)-2\tss
t_{31}^{(1)}\tss T_{11}(u)
-t_{21}^{(1)}\tss T_{21}(u).
\een
Together with the relation
\beql{tooxi}
T_{11}(u)\ts\xi=(u+\al-1/2)(u+\al)\ts\xi
\eeq
they imply the claim for the operator $T_{11}(u)$. For the remaining operators the
property follows from the relations
\ben
[t_{12}^{(1)},T_{21}(u)]=T_{11}(u)-T_{22}(u),\qquad
[t_{21}^{(1)},T_{22}(u)]=T_{32}(u)-T_{21}(u)
\een
and
\ben
[t_{12}^{(1)},T_{11}(u)]=T_{12}(u),\qquad
[t_{23}^{(1)},T_{32}(u)]=T_{33}(u)-T_{22}(u),\qquad
[t_{23}^{(1)},T_{33}(u)]=-T_{23}(u),
\een
which are consequences of \eqref{defrel}.
\epf

For any $r,s\in\ZZ_+$ introduce vectors of the small Verma module $M(\al)$ by setting
\begin{multline}
\xi_{rs}=T_{21}(-\al-r+3/2)\dots T_{21}(-\al-1/2)\tss T_{21}(-\al+1/2)\\[0.3em]
{}\times T_{21}(-\al-s+1)\dots T_{21}(-\al-1)\tss T_{21}(-\al)\ts\xi.
\non
\end{multline}
We would like to show that under certain
additional conditions the vectors $\xi_{rs}$ form a basis of $M(\al)$; see Theorem~\ref{thm:rvm}
and Corollary~\ref{cor:arba} below.
This will require a few lemmas where the action of the operators $T_{ij}(u)$
on these vectors is calculated.

\ble\label{lem:toneone}
In the module $M(\al)$ we have
\ben
T_{11}(u)\ts\xi_{rs}=(u+\al+r-1/2)(u+\al+s)\ts\xi_{rs}.
\een
\ele

\bpf
The formula holds for $\xi_{00}=\xi$
by \eqref{tooxi}.
The defining relations \eqref{defrel} give
\ben
T_{11}(u)\tss T_{21}(v)=\frac{u-v+1}{u-v}\ts T_{21}(v)\tss T_{11}(u)
-\frac{1}{u-v}\ts T_{21}(u)\tss T_{11}(v),
\een
which implies the desired formula by an obvious induction.
\epf

\ble\label{lem:ttwoone}
In the module $M(\al)$ for all $r\leqslant s+1$ we have
\ben
T_{21}(u)\ts\xi_{rs}=\frac{(-1)^{r+1}\tss (s-r+1)(2u+2\al+2r-1)}{(s+1)(2s-2r+1)}\ts\xi_{r,s+1}
+\frac{2\tss(u+\al+s)}{2s-2r+1}\ts\xi_{r+1,s}.
\een
\ele

\bpf
By the definition of the vectors $\xi_{rs}$ we have $T_{21}(-\al-r+1/2)\tss \xi_{rs}=\xi_{r+1,s}$.
Next we point out the following relation for generators of $\X(\osp_{1|2})$:
\ben
(u-v-1/2)\ts t_{21}(u)\tss t_{21}(v)+(u-v+1/2)\ts t_{21}(v)\tss t_{21}(u)
=t_{31}(v)\tss t_{11}(u)-t_{31}(u)\tss t_{11}(v).
\een
It is derived by calculating the commutators $[t_{21}(u), t_{21}(v)]$ and
$[t_{11}(u), t_{31}(v)]$ by \eqref{defrel} and eliminating the term $t_{11}(u)\tss t_{31}(v)$.
By Lemma~\ref{lem:toneone} we have $T_{11}(u)\tss \xi_{rs}=0$ for $u=-\al-r+1/2$ and
$u=-\al-s$. Hence, we come to the relation
\ben
(r-s-1)\ts T_{21}(-\al-s)\tss T_{21}(-\al-r+1/2)\ts\xi_{rs}=
-(r-s)\ts T_{21}(-\al-r+1/2)\tss T_{21}(-\al-s)\ts\xi_{rs}.
\een
Since $T_{21}(-\al-s)\ts\xi_{0\tss s}=\xi_{0,s+1}$,
applying the relation repeatedly, we get the formula
\beql{rls}
T_{21}(-\al-s)\ts\xi_{rs}=\frac{(-1)^r\tss (s-r+1)}{s+1}\ts \xi_{r,s+1}
\eeq
which is valid for all $r\leqslant s+1$. Finally, using
the Lagrange interpolation formula
\ben
T_{21}(u)=\frac{u+\al+r-1/2}{r-s-1/2}\ts T_{21}(-\al-s)-\frac{u+\al+s}{r-s-1/2}
\ts T_{21}(-\al-r+1/2).
\een
we get the relation in the lemma.
\epf

\ble\label{lem:tonetwo}
In the module $M(\al)$ for all $r\leqslant s$ we have
\ben
\bal
T_{12}(u)\ts\xi_{rs}&=-\frac{r\tss (s-r+1)(2\tss\al+2r-3)(u+\al+s)}{2\tss(2s-2r+1)}\ts\xi_{r-1,s}
\\[0.4em]
{}&+\frac{(-1)^{r+1}\tss s\tss (2s+1)(\al+s-1)(2u+2\al+2r-1)}{4\tss(2s-2r+1)}\ts\xi_{r,s-1}.
\eal
\een
\ele

\bpf
By Proposition~\ref{prop:pol}, the operator $T_{12}(u)$ is a polynomial in $u$ of degree one.
As in the proof of Lemma~\ref{lem:ttwoone}, it will be sufficient to calculate the action
of the operator for two different values $u=-\al-r+1/2$ and $u=-\al-s$,
and then apply the Lagrange interpolation
formula.

Recall from Sec.~\ref{sec:gd} that the coefficients of the series $h_1(u)$ and $h_2(u)$
pairwise commute. Set $d(u)=h_1(u)\tss h_2(u+1)$.
Using the defining relations \eqref{defrel}, we
can also write this series in the form
\ben
d(u)=t_{22}(u)\tss t_{11}(u+1)+t_{12}(u)\tss t_{21}(u+1).
\een
The coefficients of the series $c(u)$ act by scalar multiplication in the small Verma module.
The scalars are found from \eqref{cu} and given by
\beql{cuac}
c(u)\mapsto \frac{(u+1)\tss (u+\al)}{u\tss(u+\al+1)}.
\eeq
On the other hand, by Lemma~\ref{lem:toneone},
the coefficients of the series $h_1(u)=t_{11}(u)$ act on each vector $\xi_{rs}$
as multiplications by scalars depending on $r$ and $s$.
Hence the same property holds for the coefficients
of $d(u)$ whose action is uniquely determined by the relation
\ben
d(u)\tss d(u+1/2)=c(u)\tss h_1(u+1/2)\tss h_1(u+1)
\een
implied by
\eqref{cu}. Therefore, the action is found by
\ben
d(u)\mapsto \frac{(u+1/2)\tss (u+\al)}{u\tss(u+\al+1/2)}\ts h_1(u+1/2).
\een
For the corresponding polynomial operator
\beql{Du}
D(u)=T_{22}(u)\tss T_{11}(u+1)+T_{12}(u)\tss T_{21}(u+1)
\eeq
we then have
\beql{Dud}
D(u)=(u+1)(u+\al-1/2)\ts T_{11}(u+1/2).
\eeq
For any $r,s\in\ZZ_+$ we find from \eqref{Du} by applying Lemma~\ref{lem:toneone} that
\ben
\bal
D(-\al-r-1/2)\ts\xi_{rs}&=T_{12}(-\al-r-1/2)\tss T_{21}(-\al-r-1/2)\ts \xi_{rs}\\
{}&=T_{12}(-\al-r-1/2)\ts \xi_{r+1,s}.
\eal
\een
Hence using \eqref{Dud} and replacing $r$ by $r-1$ we find
\ben
T_{12}(-\al-r+1/2)\ts \xi_{rs}=-\frac{1}{4}\ts r(s-r+1)(2\tss\al+2r-3)\ts\xi_{r-1,s},
\een
which holds for $r\geqslant 1$. To extend this formula to the case $r=0$ use Lemma~\ref{lem:toneone}
and relations
\beql{toott}
[T_{12}(u)\tss T_{21}(v)]=\frac{1}{u-v}\ts \big(T_{22}(u)\tss T_{11}(v)-T_{22}(v)\tss T_{11}(u)\big)
\eeq
implied by \eqref{defrel} to derive by induction on $s$ that
$T_{12}(-\al+1/2)\ts \xi_{0\tss s}=0$.

Similarly, taking $u=-\al-s-1$ in \eqref{Du} and \eqref{Dud}, we get by using \eqref{rls} that
\ben
T_{12}(-\al-s)\ts \xi_{rs}=\frac{1}{4}\ts (-1)^r\tss s\tss(2s+1)(\al+s-1)\ts\xi_{r,s-1},
\een
which holds for $r<s$. This formula extends to the case $r=s$ by applying relation \eqref{toott}
and taking into account Lemma~\ref{lem:toneone}.
\epf

\bth\label{thm:rvm}
Suppose that $-\al\notin \ZZ_+$ and $-\al+1/2\notin\ZZ_+$. Then the $\X(\osp_{1|2})$-module $M(\al)$
is irreducible. Moreover, the vectors $\xi_{rs}$ with $r\leqslant s$ form
a basis of $M(\al)$ and $\xi_{rs}=0$ for $r>s$.
\eth

\bpf
We start by showing that all vectors $\xi_{rs}$ with $0\leqslant r\leqslant s$ are nonzero
in $M(\al)$. The conditions on $\al$ and Lemma~\ref{lem:tonetwo} imply that
it is sufficient to verify that $\xi\ne 0$; the vector $\xi_{rs}$ would then also have to be
nonzero, because the application of suitable operators $T_{12}(v)$ to $\xi_{rs}$
gives the vector $\xi$ with a nonzero coefficient.

The relation $\xi=0$ in $M(\al)$ would mean
that $\xi$, as an element of the Verma module $M(\la(u))$ with the highest weight
given in \eqref{hwal}, belongs to the submodule $K$. That is, $\xi$ is a linear
combination of vectors of the form
\ben
x_r\tss t_{21}^{(r)}\tss\xi\quad\text{for}\quad r\geqslant 2\Fand
y_r\tss\big(t_{31}^{(r)}+(\al-1/2)\tss t_{31}^{(r-1)}\big)\tss\xi\quad\text{for}\quad r\geqslant 3,
\een
with $x_r,y_r\in\X(\osp_{1|2})$. The elements $x_r$ and $y_r$ must have the respective
$\osp_{1|2}$-weights $1$ and $2$ as eigenvectors of the operator $F_{11}$.
Write these elements as linear combinations of the vectors of
the Poincar\'e--Birkhoff--Witt basis of $\X(\osp_{1|2})$
by using any ordering on the generators consistent with the increasing $\osp_{1|2}$-weights.
The right-most generators occurring in each basis monomial will have
positive $\osp_{1|2}$-weights.
On the other hand, calculating in the Verma module $M(\la(u))$ we find
\ben
\bal
t_{12}(u)\tss \big(t_{21}(v)-v^{-1}t_{21}^{(1)}\big)\ts\xi{}&=
\frac{1}{u-v}\ts \big(t_{22}(u)\tss t_{11}(v)-t_{22}(v)\tss t_{11}(u)\big)\ts\xi\\[0.3em]
{}&-v^{-1}\ts \big(t_{11}(u)-t_{22}(u)\big)\ts\xi=0,
\eal
\een
as the coefficient of $\xi$ equals
\ben
\frac{1}{u-v}\ts \big(\frac{v+\al}{v}-\frac{u+\al}{u}\big)-\al\tss u^{-1}v^{-1}=0.
\een
Now combine the second family of generators of the submodule $K$
given in \eqref{gensub} into the generating series
\ben
t_{31}(v)-v^{-1}\tss t_{31}^{(1)}
-\frac{1}{v(v+\al-1/2)}\ts t_{31}^{(2)}
\een
which can be written as the anti-commutator of $t_{21}^{(1)}$ with the series
\ben
t_{21}(v)-v^{-1}\tss t_{21}^{(1)}
-\frac{1}{v(v+\al-1/2)}t_{21}^{(2)}
\een
whose coefficients are also generators of $K$. Working first with one part of the
anti-commutator and using the previous calculation we get
\ben
t_{12}(u)\tss t_{21}^{(1)}\tss \big(t_{21}(v)-v^{-1}\tss t_{21}^{(1)}\big)\ts\xi
=\big(t_{11}(u)-t_{22}(u)\big)\big(t_{21}(v)-v^{-1}\tss t_{21}^{(1)}\big)\ts\xi.
\een
By the previous argument, the coefficients of this series vanish under
the action of the coefficients of the series $t_{12}(w)$.
Turning to the second part of the anti-commutator, we find that the expression
\ben
t_{12}(u)\tss \big(t_{21}(v)-v^{-1}\tss t_{21}^{(1)}
-\frac{1}{v(v+\al-1/2)}\tss t_{21}^{(2)}\big)\tss t_{21}^{(1)}\ts \xi
\een
equals
\beql{paro}
-\big(t_{21}(v)-v^{-1}\tss t_{21}^{(1)}
-\frac{1}{v(v+\al-1/2)}t_{21}^{(2)}\big)\tss t_{12}(u)\tss t_{21}^{(1)}\ts \xi
\eeq
plus
\ben
\bal
\frac{1}{u-v}\ts \big(t_{22}(u)\tss t_{11}(v)&-t_{22}(v)\tss t_{11}(u)\big)\tss t_{21}^{(1)}\ts\xi
-v^{-1}\ts \big(t_{11}(u)-t_{22}(u)\big)\tss t_{21}^{(1)}\ts\xi\\
{}&-\frac{1}{v(v+\al-1/2)}\tss \big((u+t_{22}^{(1)})\tss t_{11}(u)-t_{22}(u)(u+t_{11}^{(1)})\big)
\tss t_{21}^{(1)}\ts\xi.
\eal
\een
The expression \eqref{paro} vanishes under
the action of the coefficients of the series $t_{12}(w)$, so we only need to transform
the second expression. We will do this modulo terms of the form $x_r t_{21}^{(r)}\tss\xi$
with $r\geqslant 2$ which were already considered above.
Note the commutators
\ben
[t_{11}(u),t_{21}^{(1)}]=t_{21}(u),\qquad [t_{22}(u),t_{21}^{(1)}]=t_{21}(u)-t_{32}(u).
\een
Using the second relation in \eqref{ef} and writing the Gaussian
generators in terms of the $t_{ij}(u)$, we find
\ben
t_{21}(u)\tss t_{22}(u+1/2)\ts\xi=t_{32}(u)\tss t_{11}(u+1/2)\ts\xi.
\een
Since $t_{21}(u)\tss\xi\equiv u^{-1}\tss t_{21}^{(1)}\ts\xi$, we derive that
$t_{32}(u)\tss\xi\equiv (u+\al-1/2)^{-1}\tss t_{21}^{(1)}\ts\xi$. Therefore,
the expression in question is then simplified
by using relations
\ben
t_{11}(u)\tss t_{21}^{(1)}\ts\xi\equiv u^{-1}\tss t_{21}^{(1)}\ts\xi\Fand
t_{22}(u)\tss t_{21}^{(1)}\ts\xi\equiv \frac{u^2+(\al-1/2)(u+1)}{u(u+\al-1/2)}\tss t_{21}^{(1)}\ts\xi
\een
and thus verifying that it reduces to zero. This completes the proof that $\xi\not\equiv 0\mod K$.

As a next step, observe that since the vectors $\xi_{rs}$ with $0\leqslant r\leqslant s$ are nonzero
in $M(\al)$, they are eigenvectors for the operator $T_{11}(u)$, whose eigenvalues
are distinct as polynomials in $u$. Hence they are linearly independent. The number
of those vectors of the $\osp_{1|2}$-weight $-\al-p$ equals $\lfloor p/2\rfloor+1$,
which together with the inequality \eqref{weineq} proves that they form a basis
of the weight space $M(\al)_{-\al-p}$. Thus, all vectors
$\xi_{rs}$ with $0\leqslant r\leqslant s$ form a basis of $M(\al)$.
Any vector $\xi_{rs}$ with $r>s$ cannot be nonzero, because otherwise it would be an eigenvector
for the operator $T_{11}(u)$ whose eigenvalue does not occur among those of the vectors in $M(\al)$.

Finally, we prove the irreducibility of $M(\al)$. As we noted
in the beginning of the proof,
the application of suitable operators $T_{12}(v)$ to an arbitrary basis vector $\xi_{rs}$
yields the highest vector $\xi$ with a nonzero coefficient. This implies that any nonzero
submodule of $M(\al)$ must contain $\xi$ and so coincide with $M(\al)$.
\epf

\bco\label{cor:arba}
For any $\al\in\CC$ the vectors $\xi_{rs}$ with $0\leqslant r\leqslant s$ form
a basis of $M(\al)$.
\eco

\bpf
Consider the vector space $\wt M(\al)$ with basis elements $\wt\xi_{rs}$ labelled
by $r,s\in\ZZ_+$ with $0\leqslant r\leqslant s$. Define the action
of the generators $t_{11}^{(r)}$, $t_{21}^{(r)}$ and $t_{12}^{(r)}$ of $X(\osp_{1|2})$
in $\wt M(\al)$ by using the formulas of Lemmas~\ref{lem:toneone}, \ref{lem:ttwoone} and
\ref{lem:tonetwo}, where the vectors $\xi_{rs}$ with $r\leqslant s$
are respectively replaced with $\wt\xi_{rs}$,
while all vectors $\xi_{rs}$ with $r>s$ are replaced by $0$.
Also, let the coefficients of the series $c(u)$ act in $\wt M(\al)$ by scalar multiplication
defined by \eqref{cuac}. By Theorem~\ref{thm:rvm},
this assignment endows the space $\wt M(\al)$ with a $X(\osp_{1|2})$-module structure
for all $-\al\notin \ZZ_+$ and $-\al+1/2\notin\ZZ_+$. Since the matrix elements of the generators
in the basis depend polynomially on $\al$, the same formulas define a representation of
$X(\osp_{1|2})$ in $\wt M(\al)$ for all values of $\al$ by continuity.

The formulas for the action of the generators in the basis $\wt\xi_{rs}$
show that for any $\al\in\CC$
there is an $X(\osp_{1|2})$-module epimorphism $\pi:M(\la(u))\to \wt M(\al)$
defined by $\xi\mapsto \wt\xi_{00}$,
where the highest weight $\la(u)$ of the Verma module
is given by \eqref{hwal}. Moreover, the submodule $K$ of $M(\la(u))$ is contained
in the kernel of $\pi$ which gives rise to an epimorphism $\bar\pi:M(\al)\to \wt M(\al)$
with $\xi_{rs}\mapsto \wt\xi_{rs}$.
By taking into account the dimensions of the respective $\osp_{1|2}$-weight components,
we conclude from \eqref{weineq} that $\bar\pi$ is an isomorphism.
\epf

As was pointed out in the proof of Corollary~\ref{cor:arba}, for any $\al\in\CC$
the vectors \eqref{vects}
form a basis of $M(\al)$, and \eqref{weineq} is in fact an equality:
\ben
\dim M(\al)_{-\al-p}=\lfloor p/2\rfloor+1.
\een

\subsection{Elementary modules}
\label{subsec:em}

The {\em elementary modules} $L(\al)$ are defined as the irreducible quotients of $M(\al)$.
We would like to describe
the structure of $L(\al)$ for the values of $\al$ which do not satisfy the
assumptions of Theorem~\ref{thm:rvm}; that is,
$-\al\in \ZZ_+$ or $-\al+1/2\in\ZZ_+$.

\bpr\label{prop:fin}
Suppose that $-\al=k\in \ZZ_+$. The linear span $J$ of all basis vectors $\xi_{rs}$ of $M(-k)$
with $s>k$ is an $\X(\osp_{1|2})$-submodule.
The module $L(-k)$ is isomorphic to the quotient $M(-k)/J$, and
the vectors $\xi_{rs}\mod J$ with $0\leqslant r\leqslant s\leqslant k$ form
its basis.
\epr

\bpf
The formula of Lemma~\ref{lem:tonetwo} gives
\ben
T_{12}(u)\ts\xi_{r,k+1}=\frac12\ts r\tss (k-r+2)(u+1)\ts\xi_{r-1,k+1}
\een
for all $r\leqslant k+1$. This implies that the subspace $J$ of $M(-k)$
is invariant under the action of $\X(\osp_{1|2})$. Furthermore,
the formula of Lemma~\ref{lem:tonetwo} also shows that the quotient $M(-k)/J$
is irreducible and hence isomorphic to $L(-k)$.
\epf

\bpr\label{prop:infin}
Suppose that $-\al+1/2=k\in \ZZ_+$. The linear span
$I$ of all basis vectors $\xi_{rs}$ of $M(-k+1/2)$
with $r>k$ is an $\X(\osp_{1|2})$-submodule.
The module $L(-k+1/2)$ is isomorphic to the quotient $M(-k+1/2)/I$, and
the vectors $\xi_{rs}\mod I$ with $0\leqslant r\leqslant k$ form
its basis.
\epr

\bpf
The formula of Lemma~\ref{lem:tonetwo} now gives
\ben
T_{12}(u)\ts\xi_{k+1,s}
=\frac{(-1)^{k}}{4}\tss s\tss (2s+1)(u+1)\ts\xi_{k+1,s-1}
\een
for all $s\geqslant k+1$. Recalling that $\xi_{rs}=0$ for $r>s$ we conclude that
the subspace $I$ of $M(-k+1/2)$
is invariant under the action of $\X(\osp_{1|2})$. Furthermore,
Lemma~\ref{lem:tonetwo} implies that the quotient $M(-k+1/2)/I$
is irreducible and hence isomorphic to $L(-k+1/2)$.
\epf

\bco\label{cor:irrfin}
We have the following criteria.
\begin{enumerate}
\item
The $\X(\osp_{1|2})$-module $M(\al)$ is irreducible if and only if
$-\al\notin \ZZ_+$ and $-\al+1/2\notin\ZZ_+$.
\item
The $\X(\osp_{1|2})$-module $L(\al)$ is finite-dimensional if and only if
$-\al=k\in \ZZ_+$. Moreover,
\ben
\dim L(-k)=\binom{k+2}{2}.
\een
\end{enumerate}
\eco

\bpf
All parts are immediate from Theorem~\ref{thm:rvm} and Propositions~\ref{prop:fin} and \ref{prop:infin}.
\epf

As the above description of the elementary modules shows, they admit bases
formed by $\osp_{1|2}$-weight vectors. Accordingly, we can define their {\em characters}
by using formal exponents of a variable $q$ and using the definition \eqref{wei}
of
$\osp_{1|2}$-weight subspaces. Namely, we set
\ben
\ch V=\sum_{\ga} \dim V_{-\ga}\ts q^{\ga}.
\een
In particular, the character of the
irreducible highest weight module $V(\mu)$ over
$\osp_{1|2}$ is found by
\ben
\ch V(\mu)=\frac{\ \ q^{-\mu}}{1-q}\Fand\ch V(\mu)=\frac{q^{-\mu}-q^{\tss\mu+1}}{1-q}
\een
for $\mu\notin\ZZ_+$ and $\mu\in\ZZ_+$, respectively.

\bco\label{cor:char}
\begin{enumerate}
\item
The character of $M(\al)$ is given by
\ben
\ch M(\al)=\frac{q^{\al}}{(1-q)(1-q^2)}.
\een
\item
For $-\al=k\in \ZZ_+$ we have
\ben
\ch L(-k)=q^{-k}\ts\frac{(1-q^{k+1})(1-q^{k+2})}{(1-q)(1-q^2)}.
\een
\item
For $-\al+1/2=k\in \ZZ_+$ we have
\ben
\ch L(-k+1/2)=q^{-k+1/2}\ts\frac{1-q^{2k+2}}{(1-q)(1-q^2)}.
\een
\end{enumerate}
\eco

\bpf
The formulas follow by evaluating the dimensions of the weight subspaces.
\epf

In terms of the characters of the $\osp_{1|2}$-modules, we can write
the above formulas as
\ben
\ch L(-k)=\sum_{p=0}^{\lfloor k/2\rfloor} \ch V(k-2p)
\een
and
\ben
\ch L(-k+1/2)=\sum_{p=0}^{k} \ch V(k-1/2-2p).
\een

Finite-dimensional modules over the Lie superalgebras $\osp_{1|2n}$
are known to be completely reducible;
see e.g. \cite[Sec.~2.2.5]{cw:dr}.
The formulas for the action of the generator $F_{12}$ of $\osp_{1|2}$ in the basis
$\xi_{rs}$ of $L(-k)$ show that there are singular vectors of the weights $k$, $k-2$, etc.,
to imply the direct sum decomposition
\ben
L(-k)\cong\bigoplus_{p=0}^{\lfloor k/2\rfloor} V(k-2\tss p).
\een

\bco\label{cor:res}
The restriction of the module $L(\al)$ to the Lie superalgebra $\osp_{1|2}$
is irreducible if and only if $\al=0,-1$ or $1/2$.
\qed
\eco

Corollary~\ref{cor:res} shows that the $\osp_{1|2}$-modules $V(0)$, $V(1)$
and $V(-1/2)$ can be extended to $\X(\osp_{1|2})$. The Yangian action
on the three-dimensional vector representation $V(1)=\CC^{1|2}$ arises from
the replacement of $T(u)$ in
the $RTT$-relation \eqref{RTT} by a transposed $R$-matrix $R(u)$ which satisfies
the Yang--Baxter equation; cf. \cite{aacfr:sy}. It is given explicitly by setting
\beql{vectre}
t_{ij}(u)\mapsto \de_{ij}+u^{-1}\tss e_{ij}(-1)^{\bi}-(u+\ka)^{-1}\tss e_{j'i'}(-1)^{\bi\bj}\ts\ta_i\ta_j
\eeq
and is isomorphic to $L(-1)$.

\subsection{Tensor product modules}
\label{subsec:tp}

We will now use the results of the previous
sections to complete the proof of the Main Theorem.

Recall that the elementary modules of the form $L(\al,\be)$
and small Verma modules $M(\al,\be)$ are
associated with the highest weights of the form \eqref{hwelem}.
They can be obtained by twisting the respective modules $L(\al)$
and $M(\al)$ with the shift automorphisms \eqref{shift}.
Corollary~\ref{cor:irrfin}\tss(2) implies that the module $L(\al,\be)$
is finite-dimensional if and only if $\be-\al\in \ZZ_+$.

For the highest weight of the form \eqref{hwgen},
the existence of a monic polynomial $P(u)$ satisfying \eqref{ydomire}
is equivalent to the condition that the parameters $\be_1,\dots,\be_k$
can be renumbered in such a way that all differences $\be_i-\al_i$ with
$i=1,\dots,k$ belong to $\ZZ_+$. If this condition holds, then
the tensor product module \eqref{tenprelem} is finite-dimensional
and so is its irreducible subquotient $L(\la(u))$. This thus proves that the
conditions of the Main Theorem are sufficient for the irreducible
highest weight module to be finite-dimensional. In the rest of this section,
we will show that the conditions are also necessary.

By the results of Sec.~\ref{subsec:rv}, each small Verma
module $M(\al,\be)$ has the basis $\xi_{rs}$ parameterized
by $r,s\in\ZZ_+$ with $r\leqslant s$ and the generators of the extended Yangian $\X(\osp_{1|2})$
act by the rules implied by Lemmas~\ref{lem:toneone}, \ref{lem:ttwoone} and \ref{lem:tonetwo}.
Namely, for all $i,j\in\{1,2,3\}$ we now introduce the operators
$T_{ij}(u)=(u+\al-1/2)\tss (u+\be)\ts t_{ij}(u)$, and
the formulas take the following form, where the vectors $\xi_{rs}$ with $r>s$
are equal to zero:
\ben
T_{11}(u)\ts\xi_{rs}=(u+\al+r-1/2)(u+\al+s)\ts\xi_{rs}
\een
together with
\ben
T_{21}(u)\ts\xi_{rs}=\frac{(-1)^{r+1}\tss (s-r+1)(2u+2\al+2r-1)}{(s+1)(2s-2r+1)}\ts\xi_{r,s+1}
+\frac{2\tss(u+\al+s)}{2s-2r+1}\ts\xi_{r+1,s}
\een
and
\ben
\bal
T_{12}(u)\ts\xi_{rs}&=-\frac{r\tss (s-r+1)(2\tss\al-2\tss\be+2r-3)(u+\al+s)}{2\tss(2s-2r+1)}
\ts\xi_{r-1,s}\\[0.4em]
{}&+\frac{(-1)^{r+1}\tss s\tss (2s+1)(\al-\be+s-1)(2u+2\al+2r-1)}{4\tss(2s-2r+1)}\ts\xi_{r,s-1}.
\eal
\een
The coefficients of the series $c(u)$ act on $M(\al,\be)$ by scalar multiplication,
with the scalars found from \eqref{cu} and given by
\ben
c(u)\mapsto \frac{(u+\al)\tss (u+\be+1)}{(u+\al+1)\tss (u+\be)}.
\een

By Corollary~\ref{cor:irrfin}\tss(1),
the $\X(\osp_{1|2})$-module $M(\al,\be)$ is irreducible if and only if
$\be-\al\notin \ZZ_+$ and $\be-\al+1/2\notin\ZZ_+$. In the cases where
$M(\al,\be)$ is reducible, the above formulas for the action of $T_{ij}(u)$
extend to the irreducible quotients $L(\al,\be)$ with the assumption that the vectors
$\xi_{rs}$ belonging to the maximal
proper submodule of $M(\al,\be)$ are understood as equal to zero.

Our argument will rely on certain sufficient conditions for the tensor
product of the form \eqref{tenprelem} to be irreducible as
an $\X(\osp_{1|2})$-module. To state the conditions we will use a notation
involving multisets of complex numbers $\{z_1,\dots,z_l\}$. For such a multiset we will write
$\{z_1,\dots,z_l\}_+$ to denote the multiset formed by all elements $z_i$ which belong to $\ZZ_+$.

\bth\label{thm:suffc}
Suppose that for each $h=1,\dots,k-1$ the following holds:
\begin{enumerate}
\item
If the multiset $\{\be_h-\al_i,\ \be_i-\al_h\ |\ i=h,\dots,k\}_+$ is not empty,
then $\be_h-\al_h$ is a minimal
element of the multiset
$\{\be_h-\al_i,\ \be_i-\al_h,\ \be_h-\al_i+1/2,\ \be_i-\al_h+1/2\ |\ i=h,\dots,k\}_+$.
\item
If the multiset $\{\be_h-\al_i,\ \be_i-\al_h\ |\ i=h,\dots,k\}_+$ is empty and the multiset\newline
$\{\be_h-\al_i+1/2,\ \be_i-\al_h+1/2\ |\ i=h,\dots,k\}_+$ is not empty,
then $\be_h-\al_h+1/2$ is a minimal
element of this multiset.
\end{enumerate}
Then the $\X(\osp_{1|2})$-module $L$ defined in \eqref{tenprelem} is irreducible.
\eth

\bpf
We let $\xi_{rs}^{(l)}$ denote the basis vectors of the module $L(\al_l,\be_l)$
with the highest vector $\xi^{(l)}$.
Proposition~\ref{prop:pol} implies that all operators
\ben
T_{ij}(u)=\prod_{l=1}^k(u+\al_l-1/2)\tss (u+\be_l)\ts t_{ij}(u)
\een
acting in the module $L$
are polynomials in $u$.

As a first step, we will show by induction on $k$ that any vector $\ze\in L$
satisfying the condition $T_{12}(u)\tss\ze=0$ is proportional to $\xi^{(1)}\ot \dots\ot\xi^{(k)}$.
The case $k=1$ is clear so we suppose that $k\geqslant 2$. We may assume that
such a vector $\ze$ is an $\osp_{1|2}$-weight vector and write
\ben
\ze=\sum_{r,s} \xi_{rs}^{(1)}\ot \ze_{rs},\qquad \ze_{rs}\in L(\al_2,\be_2)\ot\dots\ot L(\al_k,\be_k).
\een
The sum is finite and taken over the pairs $r\leqslant s$ with
the condition that the $\xi_{rs}^{(1)}$ are basis vectors
of $L(\al_1,\be_1)$.
Let $p$ be the maximal sum $r+s$ for which there are nonzero elements
$\ze_{rs}$ in the expression. By taking the coefficient of $\xi^{(1)}_{rs}$
with $r+s=p$ in the relation $T_{12}(u)\tss\ze=0$, we get $T_{12}(u)\tss \ze_{rs}=0$.
By the induction hypothesis, $\ze_{rs}$
is proportional to the vector $\xi'=\xi^{(2)}\ot \dots\ot\xi^{(k)}$. Furthermore,
the defining relations \eqref{defrel} give
\ben
T_{12}(u)\tss T_{11}(v)=\frac{u-v-1}{u-v}\ts T_{11}(v)\tss T_{12}(u)
+\frac{1}{u-v}\ts T_{11}(u)\tss T_{12}(v).
\een
Hence, for any value of $v$, the vector $T_{11}(v)\tss\ze$ is also annihilated by
the operator $T_{12}(u)$.
Note that the basis vectors $\xi_{rs}^{(1)}$ are eigenvectors
for the operator $T_{11}(v)$ with distinct eigenvalues, as polynomials in $v$.
This implies that by taking a suitable value of $v$, we can find a
linear combination of the vectors $T_{11}(v)^m\tss\ze$
with $m=0,1,\dots$ to get
an $\osp_{1|2}$-weight vector $\ze$ of the form
\beql{sive}
\ze=\xi_{r_0s_0}^{(1)}\ot\xi'+\sum_{r+s<p} \xi_{rs}^{(1)}\ot \ze_{rs},
\eeq
with $r_0+s_0=p$ such that $T_{12}(u)\tss\ze=0$.

Next we will show that the condition $r_0<s_0$ is impossible in such a vector.
Indeed, if this condition holds, consider the coefficient of
the vector $\xi_{r_0,s_0-1}^{(1)}\ot\xi'$ in the relation $T_{12}(u)\tss\ze=0$.
This coefficient can only arise from the terms
\ben
T_{12}(u)\ts \xi_{r_0,s_0}^{(1)}\ot T_{22}(u)\ts \xi'\pm
T_{11}(u)\ts \xi_{r_0,s_0-1}^{(1)}\ot T_{12}(u)\ts \ze_{r_0,s_0-1}
\een
with the sign depending on the parity of the vector $\xi_{r_0,s_0-1}^{(1)}$.
The $\osp_{1|2}$-weight condition implies that
\ben
\ze_{r_0,s_0-1}=\sum_{l=2}^k c_l\ts \xi^{(2)}\ot\dots\ot \xi_{01}^{(l)}
\ot\dots\ot \xi^{(k)}
\een
for some constants $c_l\in \CC$. We have
\ben
T_{12}(u)\ts\ze_{r_0,s_0-1}=\sum_{l=2}^k \pm\tss c_l\ts T_{11}(u)\ts\xi^{(2)}\ot\dots\ot
T_{12}(u)\ts\xi_{01}^{(l)}
\ot\dots\ot T_{22}(u)\ts\xi^{(k)}.
\een
By using the formulas for the action
of the operators $T_{ij}(u)$ and
equating the coefficient in question to zero, we get
\begin{multline}
b\tss (u+\al_1+r_0-1/2)\tss\prod_{i=2}^k(u+\al_i-1/2)(u+\be_i)\\
+(u+\al_1+r_0-1/2)(u+\al_1+s_0-1)
\sum_{l=2}^k b_l\ts\prod_{i=2}^{l-1}(u+\al_i-1/2)(u+\al_i)\\
{}\times(u+\al_l-1/2)
\prod_{i=l+1}^{k}(u+\al_i-1/2)(u+\be_i)=0,
\non
\end{multline}
where $b_l$ are some constants, while $b$ is a nonzero constant because of the condition
$s_0\leqslant \be_1-\al_1$ in the case $\be_1-\al_1\in\ZZ_+$ implied by
Proposition~\ref{prop:fin}.
By cancelling the common factors and setting $u=-\al_1-s_0+1$ we get
\ben
\prod_{i=2}^k(\be_i-\al_1-s_0+1)=0.
\een
It follows from this relation that the multiset
$\{\be_i-\al_1\ |\ i=1,\dots, k\}_+$ is not empty,
because $\be_i-\al_1=s_0-1\in\ZZ_+$ for some $i\in\{2,\dots,k\}$.
By assumption (1) of the theorem, we have
$\be_1-\al_1\in\ZZ_+$ and $\be_1-\al_1\leqslant \be_i-\al_1$.
However, this makes a contradiction, as
by Proposition~\ref{prop:fin} we must have $s_0\leqslant \be_1-\al_1$.

Excluding the condition $r_0<s_0$ in \eqref{sive}, we show next that
the condition $r_0=s_0\geqslant 1$ is impossible either.
If this condition holds, consider the coefficient of
the vector $\xi_{r_0-1,r_0}^{(1)}\ot\xi'$ in the relation $T_{12}(u)\tss\ze=0$.
This coefficient can only arise from the terms
\ben
T_{12}(u)\ts \xi_{r_0,r_0}^{(1)}\ot T_{22}(u)\ts \xi'\pm
T_{11}(u)\ts \xi_{r_0-1,r_0}^{(1)}\ot T_{12}(u)\ts \ze_{r_0-1,r_0}.
\een
By the $\osp_{1|2}$-weight condition,
\ben
\ze_{r_0-1,r_0}=\sum_{l=2}^k c_l\ts \xi^{(2)}\ot\dots\ot \xi_{01}^{(l)}
\ot\dots\ot \xi^{(k)}
\een
for some constants $c_l\in \CC$. Calculating as in the previous case,
we now come to the relation
\begin{multline}
b\tss (u+\al_1+r_0)\tss\prod_{i=2}^k(u+\al_i-1/2)(u+\be_i)\\
+(u+\al_1+r_0-3/2)(u+\al_1+r_0)
\sum_{l=2}^k b_l\ts\prod_{i=2}^{l-1}(u+\al_i-1/2)(u+\al_i)\\
{}\times(u+\al_l-1/2)
\prod_{i=l+1}^{k}(u+\al_i-1/2)(u+\be_i)=0,
\non
\end{multline}
where $b_l$ are some constants, while $b$ is a nonzero constant. The latter property
holds
because of the condition
$r_0\leqslant \be_1-\al_1+1/2$ in the case $\be_1-\al_1+1/2\in\ZZ_+$ implied by
Proposition~\ref{prop:infin}.
Cancel the common factors and set $u=-\al_1-r_0+3/2$ to get
\ben
\prod_{i=2}^k(\be_i-\al_1-r_0+3/2)=0.
\een
This means that for some $i\in\{2,\dots,k\}$ we have $\be_i-\al_1+1/2=r_0-1\in\ZZ_+$.
If the multiset $\{\be_1-\al_j,\ \be_j-\al_1\ |\ j=1,\dots,k\}_+$
is not empty,
then by assumption (1) of the theorem, we have
$\be_1-\al_1\in\ZZ_+$ and $\be_1-\al_1\leqslant \be_i-\al_1+1/2$.
This is impossible because
by Proposition~\ref{prop:fin} we must have $r_0\leqslant \be_1-\al_1$.
Hence assumption (2) of the theorem for $h=1$ should apply,
and we have $\be_1-\al_1+1/2\in\ZZ_+$ together with the inequality
\ben
\be_1-\al_1+1/2\leqslant \be_i-\al_1+1/2.
\een
This makes a contradiction, as
by Proposition~\ref{prop:infin} we must have $r_0\leqslant \be_1-\al_1+1/2$.

We have thus showed that any vector $\ze\in L$
with $T_{12}(u)\tss\ze=0$ is proportional to $\xi^{(1)}\ot \xi'$.
By looking at the set of $\osp_{1|2}$-weights
of any nonzero submodule of $L$ we derive that such a submodule
must contain a nonzero vector $\ze$ with $T_{12}(u)\tss\ze=0$,
and so contain the vector $\xi^{(1)}\ot \xi'$.
It remains to prove this vector is cyclic in $L$.

Consider the vector space $L^*$ dual to $L$ which is
spanned by all linear maps $\sigma:L\to \CC$
satisfying the condition that the linear span of the vectors $\eta\in L$
such that $\sigma(\eta)\ne 0$, is finite-dimensional.
Equip $L^*$ with an $\X(\osp_{1|2})$-module structure
by setting
\beql{dumo}
(x\ts\sigma)(\eta)=\sigma(\om(x)\ts\eta)
\quad
\text{for}
\quad
x\in\X(\osp_{1|2})
\fand
\sigma\in L^*,\ \eta\in L,
\eeq
where $\om$ is the anti-automorphism of the algebra $\X(\osp_{1|2})$
defined in \eqref{om}.
It is easy to verify that $L^*$ is isomorphic to
the tensor
product module
\beql{dul}
L(-\be_1,-\al_1)\ot \dots\ot L(-\be_k,-\al_k).
\eeq
Moreover, the highest vector of the module $L(-\be_i,-\al_i)$
can be identified with the dual basis vector $\xi^{(i)\tss *}$.
Suppose now that the submodule $N=\X(\osp_{1|2})(\xi^{(1)}\ot \dots\ot\xi^{(k)})$
of $L$ is proper and consider
its annihilator
\beql{ann}
\Ann N=\{\rho\in L^*\ |\ \rho(\eta)=0
\quad\text{for all}\quad\eta\in N\}.
\eeq
Then $\Ann N$
is a nonzero submodule of $L^*$, which does not
contain the vector $\xi^{(1)\tss *}\ot\dots\ot \xi^{(k)\tss *}$.
However, this contradicts the claim verified in the first part of the proof,
because the conditions on the parameters $\al_i$ and $\be_i$
stated in the theorem
will remain satisfied after we replace each $\al_i$ by $-\be_i$ and each
$\be_i$ by $-\al_i$.
\epf

\bpr\label{prop:shifts}
Suppose that the $\X(\osp_{1|2})$-module $L(\la(u))$ with
the highest weight \eqref{hwgen} is finite-dimensional.
Then for any nonnegative integers $l_1,\dots,l_k$ and $m_1,\dots,m_k$
the module $L(\la^+(u))$ with
the highest weight
\beql{hwgenshi}
\la^+_1(u)=\frac{(u+\al_1-l_1)\dots (u+\al_k-l_k)}{(u+\be_1+m_1)\dots (u+\be_k+m_k)}
\Fand \la^+_2(u)=1
\eeq
is also finite-dimensional.
\epr

\bpf
The highest weight module $L(\la^+(u))$ is isomorphic to an irreducible
subquotient of the finite-dimensional
module
\ben
L(\la(u))\ot L(\al_1-l_1,\al_1)\ot \dots \ot L(\al_k-l_k,\al_k)
\ot L(\be_1,\be_1+m_1) \ot \dots \ot L(\be_k,\be_k+m_k)
\een
and hence is finite-dimensional.
\epf

We now return to proving the Main Theorem. Let
the irreducible highest weight module $L(\la(u))$ with
the highest weight \eqref{hwgen} be finite-dimensional.
To argue by contradiction, suppose that it is impossible to renumber
the parameters $\be_1,\dots,\be_k$
in such a way that all differences $\be_i-\al_i$ with
$i=1,\dots,k$ belong to $\ZZ_+$. By Proposition~\ref{prop:shifts},
all modules $L(\la^+(u))$ with
the highest weight of the form
\eqref{hwgenshi} are also finite-dimensional. It is possible to
choose nonnegative integers $l_i$ and $m_i$
to ensure that the assumptions of Theorem~\ref{thm:suffc} are satisfied by
the shifted parameters $\al'_i=\al_i-l_i$ and $\be'_i=\be_i+m_i$, after a possible
renumbering. This can be done by induction, beginning
with the multiset
\ben
\{\be_1-\al_i,\ \be_i-\al_1\ |\ i=1,\dots,k\}
\een
and renumbering
the parameters $\al_i$ and $\be_i$, if necessary, to ensure that $\be_1-\al_1$
is a minimal element of the multiset
\beql{multo}
\{\be_1-\al_i,\ \be_i-\al_1\ |\ i=1,\dots,k\}_+
\eeq
if it is nonempty. Then assumption (1) of the theorem for $h=1$ is achieved by suitable shifts
$\al_i\mapsto\al_i-l_i$ and $\be_i\mapsto \be_i+m_i$ for $i=2,\dots,k$.
If the multiset \eqref{multo} is empty, then
assumption (2) for $h=1$ is achieved by a suitable renumbering
of the parameters $\al_i$ and $\be_i$. Then we continue in the same way
to consider the multisets for $h=2$, etc. As a result,
by Theorem~\ref{thm:suffc}, the module $L(\la^+(u))$ is
isomorphic
to the tensor product of the corresponding elementary modules.
Since it is finite-dimensional, all new differences $\be'_i-\al'_i$
must be nonnegative integers due to Corollary~\ref{cor:irrfin}\tss(2).

This argument implies, that
all the differences $\be_i-\al_i$ of the original parameters
may be assumed to be integers. Moreover,
we can apply some shifts as given in Proposition~\ref{prop:shifts}, to further suppose
that $\be_i-\al_i\in\ZZ_+$ for $i=1,\dots,k-1$,
while $\al_k-\be_k\in 1+\ZZ_+$, and that it is impossible to renumber
the parameters to make all the differences $\be_i-\al_i$ nonnegative integers.

Now consider all the parameters $\al_i$ and $\be_i$ which belong to the
$\ZZ$-coset in $\CC$ containing $\al_k$ and $\be_k$. Renumbering them, if necessary, suppose that
they correspond to $i=d+1,\dots,k$ for some $d\in\{0,1,\dots,k-1\}$.
After a further renumbering to satisfy the assumptions of
Theorem~\ref{thm:suffc}, we obtain
that the $\X(\osp_{1|2})$-module
\ben
L^{(2)}=L(\al_{d+1},\be_{d+1})\ot\dots\ot L(\al_{k},\be_{k})
\een
is irreducible. Similarly, by applying suitable shifts of Proposition~\ref{prop:shifts}
to the remaining parameters $\al_i,\be_i$
with $i=1,\dots,d$, and possible relabelling,
we may assume that they satisfy the assumptions of
Theorem~\ref{thm:suffc} and so
the $\X(\osp_{1|2})$-module
\ben
L^{(1)}=L(\al_1,\be_1)\ot\dots\ot L(\al_{d},\be_{d})
\een
is also irreducible. If the tensor product $L=L^{(1)}\ot L^{(2)}$ turns out to be
irreducible, then we arrive at a contradiction, because the module $L(\al_{k},\be_{k})$
is infinite-dimensional. So we will suppose that $L$ is not irreducible and denote by $\mu$
the $\osp_{1|2}$-weight of the vector $\xi^{(1)}\ot\dots\ot \xi^{(k)}$.
Consider the multiset
\beql{pze}
\{\be_i-\al_j+1/2\ |\ 1\leqslant i\leqslant d,\quad d+1\leqslant j\leqslant k\}_+
\eeq
and let $p_0$ denote its minimal element, if the multiset is nonempty, or set $p_0=+\infty$ otherwise.

\ble\label{lem:weispace}
The $\osp_{1|2}$-weight component $N_{\mu-p}$
of the cyclic span
\ben
N=\X(\osp_{1|2})(\xi^{(1)}\ot\dots\ot \xi^{(k)})
\een
coincides with
$L_{\mu-p}$ for all $0\leqslant p\leqslant 2\tss p_0$.
\ele

\bpf
Equip $L^*$ with an $\X(\osp_{1|2})$-module structure
by using \eqref{dumo}.
By considering the annihilator $\Ann N$, as defined by \eqref{ann}, it will be
sufficient to show that any vector $\ze\in L^*$ of the $\osp_{1|2}$-weight $\mu-p$
with the property $t_{12}(u)\tss\ze=0$
is proportional to the vector $\xi^{(1)\tss *}\ot\dots\ot \xi^{(k)\tss *}$.
As before, we will identify $L^*$ with the tensor product
module \eqref{dul} and denote
by $\check\xi^{(i)}$ the highest vector of the elementary module $L(-\be_i,-\al_i)$.
We will now follow the first part of the proof of Theorem~\ref{thm:suffc}
to derive by a reverse induction on $l\in\{1,\dots,k\}$, beginning with $l=k$, that any vector
\ben
\ze^{(l)}\in L(-\be_l,-\al_l)\ot \dots\ot L(-\be_k,-\al_k)
\een
of the $\osp_{1|2}$-weight $\mu-p$
with the property $t_{12}(u)\tss\ze^{(l)}=0$
is proportional to $\check\xi^{(l)}\ot\dots\ot \check\xi^{(k)}$.
This is clear for the values $l=d+1,\dots,k$, because the assumptions of
Theorem~\ref{thm:suffc} are satisfied by the corresponding parameters.

Now suppose that $l\in\{1,\dots,d\}$ and repeat the argument
of the first part of the proof of Theorem~\ref{thm:suffc} to come to the
expression
\ben
\ze^{(l)}=\check\xi_{r_0s_0}^{(l)}\ot\xi'+\sum_{r+s<p} \check\xi_{rs}^{(l)}\ot \ze_{rs},
\een
analogous to \eqref{sive}, where $r_0+s_0=p$ and $\xi'=\check\xi^{(l+1)}\ot\dots\ot \check\xi^{(k)}$.
Arguing as in that proof, we find that
the condition $r_0<s_0$ is impossible, leading to the only possibility
that $r_0=s_0\geqslant 1$. In this case, with our conditions of the parameters,
we must have $p=2r_0$ and
\beql{refbr}
\be_l-\al_j+1/2=r_0-1
\eeq
for some $d+1\leqslant j\leqslant k$. Since $r_0-1\in\ZZ_+$,
relation \eqref{refbr} implies that $p_0$ has a finite value and
$r_0>p_0$.
This makes a contradiction, because $p=2r_0\leqslant 2p_0$ by the
assumption, thus completing the proof of the lemma.
\epf

For any $s\in\ZZ_+$ set $\eta_s=\xi^{(1)}\ot\dots\ot \xi^{(k-1)}\ot \xi^{(k)}_{0\tss s}
\in L(\al_1,\be_1)\ot\dots\ot L(\al_k,\be_k)$.

\ble\label{lem:inff}
In the tensor product module, for any $s\in\ZZ_+$ we have
\begin{align}\label{ttone}
T_{21}(-\al_k-s)\tss \eta_s
{}&=\prod_{i=1}^{k-1}(\be_i-\al_k-s)(\al_i-\al_k-s-1/2)\ts
\eta_{s+1},\\[0.2em]
T_{12}(u)\tss \eta_s
&=\frac{s}{2}\ts(\be_k-\al_k-s+1)\ts\prod_{i=1}^{k-1}(u+\al_i-1/2)(u+\al_i)\ts
\eta_{s-1},
\non\\
\intertext{and}
T_{11}(u)\tss \eta_s
&=(u+\al_k-1/2)(u+\al_k+s)\ts\prod_{i=1}^{k-1}(u+\al_i-1/2)(u+\al_i)\ts
\eta_{s}.
\non
\end{align}
\ele

\bpf
All relations are immediate from the coproduct rule \eqref{Delta}
and the formulas for the action of the generators
of the extended Yangian
in the basis $\xi_{rs}$ of the elementary module $L(\al,\be)$, which were
recalled in
the beginning of this section. In particular, for \eqref{ttone} we take into account
the relations $T_{11}(-\al_k-s)\ts\xi^{(k)}_{0\tss s}=0$ and
$T_{21}(-\al_k-s)\ts\xi^{(k)}_{0\tss s}=\xi^{(k)}_{0, s+1}$ in $L(\al_k,\be_k)$.
\epf

Observe that the numerical coefficient on the right hand side of \eqref{ttone}
is nonzero for any values of $s$ outside the multisets
\beql{ms}
\{\be_i-\al_k\ |\ i=d+1,\dots,k-1\}_+\fand \{\al_i-\al_k-1/2\ |\ i=1,\dots,d\}_+.
\eeq
On the other hand, recalling that $p_0$ is the minimal element of the multiset \eqref{pze}
when it is nonempty,
note that we can use the shifts of the parameters $\al_i,\be_i$ with $i=1,\dots,d$
as in Proposition~\ref{prop:shifts} to keep the assumptions of Theorem~\ref{thm:suffc}
satisfied. The module $L^{(1)}$
with the shifted parameters remains irreducible, while we can make the value of $p_0$
arbitrarily large. It will be sufficient to make $p_0$ large enough for the elements
of both multisets in \eqref{ms} not to exceed $2p_0$, noting that
the elements of the second multiset can only decrease after
the shifts
$\al_i\mapsto \al_i-l_i$ for $i=1,\dots,d$.

The $\osp_{1|2}$-weight of the vector $\eta_s$ equals $\mu-s$, and hence, by
Lemma~\ref{lem:weispace}, all vectors $\eta_s$ with $s\leqslant 2p_0$
belong to the cyclic span $N=\X(\osp_{1|2})\tss\eta_0$. This property
extends to all values $s\in\ZZ_+$ by relation \eqref{ttone} of Lemma~\ref{lem:inff},
because the numerical coefficient of $\eta_{s+1}$ does not vanish for $s>2p_0$.
The remaining two relations of Lemma~\ref{lem:inff} imply that
the images of the vectors $\eta_s$ in the irreducible quotient $L(\la(u))$ of $N$
are linearly independent. Hence, $L(\la(u))$ is infinite-dimensional, as it contains
an infinite family of linearly independent
vectors. This contradiction completes the proof of the first part of the Main Theorem.
The second part concerning representations of the Yangian $\Y(\osp_{1|2})$
is immediate from the decomposition \eqref{tensordecom}; cf. \cite[Sec.~5.3]{amr:rp}.

\bigskip

Comparing the irreducibility conditions with those for the evaluation modules over
the Yangian $\Y(\gl_2)$ (see e.g. \cite[Sec.~3.3]{m:yc}), note that unlike that case,
it is not possible, in general, to renumber the parameters of the given highest weight \eqref{hwgen}
to satisfy the assumptions of Theorem~\ref{thm:suffc}. In fact, not every module
$L(\la(u))$ is isomorphic to a tensor product module of the form \eqref{tenprelem},
as illustrated by the following example.

\bex\label{ex:tpr}
To describe the $\X(\osp_{1|2})$-module $L(\la(u))$ with
\ben
\la_1(u)=\frac{(u-1)(u-5/2)}{u\tss(u-3/2)},\qquad \la_2(u)=1,
\een
consider the tensor product $L=L(-1,0)\ot L(-5/2,-3/2)$
of two three-dimensional modules.
Note that its parameters do not satisfy the assumptions of Theorem~\ref{thm:suffc}.
The module $L$ turns out to have a proper submodule $K$ which is generated by the vector
\ben
\ze=\xi_{11}^{(1)}\ot \xi^{(2)}+3\ts\xi_{01}^{(1)}\ot \xi_{01}^{(2)}-\xi^{(1)}\ot \xi_{11}^{(2)}.
\een
The submodule $K$ is one-dimensional, isomorphic to a highest weight module $L(\mu(u))$
with the components
\ben
\mu_1(u)=\mu_2(u)=\frac{(u-1/2)(u-5/2)}{(u-3/2)^2}.
\een
The module $L(\la(u))$ is isomorphic to the quotient $L/K$ with $\dim L(\la(u))=8$ and so does not
admit a tensor product decomposition of the form \eqref{tenprelem}.
\qed
\eex

To conclude, we note that by analysing submodules of reducible modules $M(\al,\be)$, we can obtain
explicit constructions of some modules $L(\la(u))$ beyond the elementary modules. In particular,
for any $k\in\ZZ_+$ the submodule of $M(-k)$ generated by the vector $\xi_{0,k+1}$
is isomorphic to the highest weight module $L(\la(u))$ with
\ben
\la_1(u)=\frac{u+1}{u}\Fand \la_2(u)=\frac{(u+1/2)(u-k-1)}{u(u-k-1/2)}.
\een
The vectors $\xi_{rs}$ with
$r\leqslant s$ and $s>k$ form its basis, and
the action of the generators is
described in Sec.~\ref{subsec:rv}. The character of $L(\la(u))$,
as defined in Sec.~\ref{subsec:em}, is found by
\ben
\ch L(\la(u))=\frac{q+q^2-q^{k+3}}{(1-q)(1-q^2)}.
\een

\bigskip\bigskip

\small

\noindent
School of Mathematics and Statistics\newline
University of Sydney,
NSW 2006, Australia\newline
alexander.molev@sydney.edu.au


\begin{thebibliography}{99}

\bibitem{aacfr:rp}
{D. Arnaudon, J. Avan, N. Cramp\'e, L. Frappat, E. Ragoucy},
{\it $R$-matrix presentation for super-Yangians $Y({\rm osp}(m\vert 2n))$},
{J. Math. Phys.}  {\bf 44}  (2003), 302--308.

\bibitem{aacfr:sy}
{D. Arnaudon, N. Cramp\'e, L. Frappat, E. Ragoucy},
{\it Super Yangian $\Y(osp(1|2))$ and the universal $R$-matrix
of its quantum double},
Comm. Math. Phys. {\bf 240} (2003), 31--51.

\bibitem{aacfr:ba}
{D. Arnaudon, J. Avan, N. Cramp\'e, A.~Doikou, L. Frappat, E. Ragoucy},
{\it Bethe ansatz equations and exact S matrices for the $osp(M|2n)$
 open super-spin chain}, Nuclear Phys. B {\bf 687} (2004), 257--278.

\bibitem{amr:rp}
{D. Arnaudon, A. Molev and E. Ragoucy},
{\it On the $R$-matrix realization of Yangians
and their representations},
Annales Henri Poincar\'e {\bf 7} (2006), 1269--1325.

\bibitem{cp:yr}
{V. Chari and A. Pressley},
{\it Yangians and $R$-matrices},
Enseign. Math. {\bf 36} (1990), 267--302.

\bibitem{cw:dr}
Sh.-J. Cheng and W. Wang,
{\it Dualities and representations of Lie superalgebras},
Graduate Studies in Mathematics, 144. AMS, Providence, RI, 2012.

\bibitem{d:nr}
{V. G. Drinfeld},
{\it A new realization of
Yangians and quantized affine algebras}, {Soviet Math. Dokl.}
{\bf 36} (1988), 212--216.

\bibitem{fikk:yy}
{J. Fuksa, A. P. Isaev, D. Karakhanyan, R. Kirschner},
{\it Yangians and Yang-Baxter $R$-operators
 for ortho-symplectic superalgebras}, Nuclear Phys. B {\bf 917} (2017), 44--85.

\bibitem{gr:dm}
{I. M. Gelfand and V. S. Retakh}, {\it Determinants
of matrices over noncommutative rings}, {Funct. Anal. Appl.} {\bf 25}
(1991), 91--102.

\bibitem{ikk:yb}
A. P. Isaev, D. Karakhanyan, R. Kirschner,
{\it Yang--Baxter $R$-operators for $osp$ superalgebras},
 Nuclear Phys. B {\bf 965} (2021), Paper No. 115355, 28 pp.

\bibitem{imo:nf}
A. P. Isaev, A. I. Molev and O. V. Ogievetsky,
{\it A new fusion procedure for the Brauer algebra and evaluation homomorphisms},
Int. Math. Res. Not. (2012), 2571--2606.

\bibitem{jl:ib}
N. Jing and M. Liu,
{\it Isomorphism between two realizations of the Yangian $Y(\mathfrak{so}_3)$},
J.~Phys.~A {\bf 46} (2013), 075201, 12 pp.

\bibitem{jlm:ib}
{N. Jing, M. Liu and A. Molev},
{\it Isomorphism between the $R$-matrix and Drinfeld
presentations of Yangian in types $B$, $C$ and $D$},
Comm. Math. Phys. {\bf 361} (2018), 827--872.

\bibitem{m:yc}
A. Molev,
{\it Yangians and classical Lie algebras}, Mathematical
Surveys and Monographs, 143. AMS,
Providence, RI, 2007.

\bibitem{m:ls}
I. M. Musson,
{\it Lie superalgebras and enveloping algebras},
Graduate Studies in Mathematics, 131. AMS, Providence, RI, 2012.

\bibitem{n:qb}
{M. L. Nazarov},
{\it Quantum Berezinian and the classical Capelli identity},
{Lett. Math. Phys.}
{\bf 21}
(1991),
123--131.

\bibitem{t:im}
{V. O. Tarasov},
{\it Irreducible monodromy matrices for the $R$-matrix of the
$XXZ$-model and lattice local quantum Hamiltonians}, {Theor. Math. Phys.}
{\bf 63} (1985),
440--454.

\bibitem{zz:rf}
{A. B. Zamolodchikov and Al. B. Zamolodchikov},
{\it Factorized $S$-matrices in two dimensions as the exact solutions
of certain relativistic quantum field models},
{Ann. Phys.} {\bf 120} (1979), 253--291.

\bibitem{zh:sy}
R. B. Zhang, {\it The $\gl(M|N)$ super Yangian and its finite-dimensional
 representations}, Lett. Math. Phys. {\bf 37} (1996), 419--434.

\end{thebibliography}
\end{document}